\documentclass[review,onefignum,onetabnum]{siamart190516}
\nolinenumbers

\usepackage[colorinlistoftodos,bordercolor=orange,backgroundcolor=orange!20,linecolor=orange,textsize=scriptsize]{todonotes}

\usepackage{graphicx}
\usepackage{tabularx}
\usepackage{multirow}

\usepackage{latexsym}

\usepackage{amsmath,amstext,amssymb,amsbsy}
\usepackage{mathtools}
\usepackage{bm,bbm}
\usepackage{caption}
\usepackage{subcaption}
\usepackage{xcolor}

\usepackage{cleveref}

\usepackage{algorithm}  
\usepackage{algorithmicx}  
\usepackage{algpseudocode}  
\usepackage{alphalph}

\usepackage[multiple]{footmisc}

\newtheorem{remark}{Remark}
\newtheorem{example}{Example}

\newcommand{\Cc}{\mathcal{C}}
\newcommand{\R}{\mathbb{R}}

\newcommand{\ba}{\boldsymbol{a}}
\newcommand{\bb}{\boldsymbol{b}}
\newcommand{\bx}{\boldsymbol{x}}
\newcommand{\bv}{\boldsymbol{v}}
\newcommand{\bp}{\boldsymbol{p}}
\newcommand{\by}{\boldsymbol{y}}
\newcommand{\bq}{\boldsymbol{q}}
\newcommand{\bu}{\boldsymbol{u}}
\newcommand{\bz}{\boldsymbol{z}}
\newcommand{\bw}{\boldsymbol{w}}

\newcommand{\bt}{\bar{t}}
\newcommand{\U}{\mathcal{U}}
\newcommand{\h}{\mathrm{h}}

\definecolor{mygreen}{rgb}{0.0, 0.5, 0.0}

\title{A time-dependent symplectic network for non-convex path planning problems with linear and nonlinear dynamics\thanks{Submitted to the editors \today.
\funding{This work was supported by
OSD/AFOSR MURI grant FA9550-20-1-0358.}}}

\author{Zhen Zhang\footnotemark[4]\footnotemark[2]\ \and
Chenye Wang\footnotemark[4]\footnotemark[2]\ \and
Shanqing Liu\footnotemark[4]\footnotemark[2]\ \and
J\'er\^ome Darbon\footnotemark[3]\footnotemark[2]  \and
George Em Karniadakis\footnotemark[2]\ 
}

\begin{document}
\maketitle
\renewcommand{\thefootnote}{\fnsymbol{footnote}}
\footnotetext[2]{Division of Applied Mathematics, Brown University 
(\{zhen\_zhang1, chenye\_wang,  shanqing\_liu, jerome\_darbon, george\_karniadakis\}@brown.edu).}
\footnotetext[3]{Corresponding author.}
\footnotetext[4]{These authors contributed equally to this work.}
\renewcommand{\thefootnote}{\arabic{footnote}}

\begin{abstract}
    We propose a novel neural network architecture (TSympOCNet) to address high--dimensional optimal control problems with linear and nonlinear dynamics. An important application of this method is to solve the path planning problem of multi-agent vehicles in real time. The new method extends our previous SympOCNet framework by introducing a time-dependent symplectic network into the architecture. In addition, we  propose a more general latent representation, which greatly improves model expressivity based on the universal approximation theorem. We demonstrate the efficacy of TSympOCNet in path planning problems with obstacle and collision avoidance, including systems with Newtonian dynamics and non-convex environments, up to dimension 512. Our method shows significant promise in handling efficiently both complex dynamics and constraints.
\end{abstract}



\section{Introduction}
\subsection{Motivation and Background}
Optimal control problems are encountered widely in practical scenarios, including path planning~\cite{Coupechoux2019Optimal,Hofer2016Application,Parzani2017HJB,Lee2021Hopf}, humanoid robot control~\cite{kuindersma2016optimization,Fujiwara2007optimal,fallon2015architecture,denk2001synthesis}, and robot manipulator control~\cite{Jin2018Robot,Kim2000intelligent,Lin1998optimal,Chen2017Reachability}. 
In the aforementioned practical applications, the control problems can rarely be solved analytically, thus necessitating approximation for the solution via numerical methods. 

One approach to address optimal control problems is the dynamic programming, pioneered by Richard Bellman~\cite{bellman1966dynamic}. 
In this approach, the optimal control problem is associated with a Hamilton-Jacobi-Bellman (HJB) equation and, in particular, the value function is characterized by the viscosity solution of the HJB PDE (see for instance~\cite{Bardi1997Optimal}). 
The HJB equation, formulated in the same dimension as the state space, encounters difficulties
in real-world applications where the state space has a high dimension.
Considering for instance a path planning problem for multiple drones. 
The state space is often formulated to be proportional (two or three times) to the number of drones in the problem. 
However, the dynamic programming approach becomes infeasible in dimensions higher than four in its original formulation.

Another approach is the maximum principle~\cite{boltyanski1960maximum, gamkrelidze2006discovery}, established by Pontryagin and colleagues almost concurrently with dynamic programming. 
In this approach, a necessary condition for the optimal trajectory is characterized by an ODE system, introducing the so-called costate which plays the role of a dual variable. 
Then, the ODE system is often numerically solved by a shooting method (see for instance~\cite{pesch1994practical}). 
Although this approach is less sensitive to the dimension, it depends strongly on the initial guess of the shooting method. The algorithm may either not converge at all (in most cases), or converge to a local optimum in the presence of nonconvexity and/or nonlinearity.

The development of efficient numerical methods for solving optimal control problems in high dimension remains a prominent research topic, and it is gaining increasing popularity with the rise in computational power. 
Various techniques have been devised in the last decade for this problem, 
for instance optimization based  methods~\cite{darbon2015convex,Darbon2016Algorithms,darbon2019decomposition,darbon2021hamilton,Chen2021Lax,Chen2021Hopf,yegorov2017perspectives,Lee2021Computationally,Kirchner2020HJ}, max-plus numerical methods~\cite{akian2008max, akian2023adaptive, dower2015maxconference,Fleming2000Max,McEneaney2006maxplus,McEneaney2007COD}, tensor decomposition techniques~\cite{dolgov2019tensor,horowitz2014linear,todorov2009efficient}, sparse grids~\cite{bokanowski2013adaptive,garcke2017suboptimal,kang2017mitigating}, polynomial approximation~\cite{kalise2019robust,kalise2018polynomial}, model order reduction~\cite{alla2017error,kunisch2004hjb}, optimistic planning~\cite{Bokanowski2021Optimistic}, hierarchical dynamic programming, reinforcement learning methods~\cite{Chen2015Exact,Chen2015Safe,alla2019efficient,akian2023multi,zhou2021actor}, and methods based on neural networks~\cite{bachouch2018deep,bansal2020deepreach,Djeridane2006Neural,jiang2016using,Han2018Solving,hure2018deep,lambrianides2019new,reisinger2019rectified,royo2016recursive,Sirignano2018DGM,Li2020generating,darbon2020overcoming,Darbon2021Neural,darbon2021neuralcontrol,nakamurazimmerer2021adaptive,jin2020sympnets,onken2021neural}. However, path planning problems are generally hard to solve directly using these algorithms, primarily due to state constraints imposed by complex environments.

Recent advances in neural networks have shown promise in overcoming the curse of dimensionality by encoding physical information into network architectures and loss functions \cite{raissi2019physics}. 
Leveraging this, SympOCNet \cite{meng2022sympocnet} was proposed to solve the Hamiltonian system corresponding to the high-dimensional optimal control problem. SympOCNet consists of two parts: (i) a latent representation, which is parameterized by affine maps and (ii) a coordinate transformation module, which is parameterized by symplectic networks (SympNets) \cite{jin2020sympnets}. Even though SympOCNet has been shown to scale to high dimensions, e.g., being able to plan the path for 256 agents, its formulation is mainly restricted to the linear integrator where the dynamics of the system does not depend on the state, i.e., controlling only the velocity of vehicles.

\subsection{Contribution}
In this paper, we consider a novel neural network architecture, the time dependent Symplectic optimal control neural network, called ``TSympOCNet'', to solve high-dimensional optimal control problems with more general, potential nonlinear,   dynamics. 
We consider a simpler Hamiltonian ODE system, that corresponds to a linear quadratic optimal control problem, in the latent space. 
We construct time-dependent symplectic maps, which are approximated and parameterized by SympNet, to transform the solution from latent space to the phase space. 
A physics-informed loss function is then constructed in the phase space, which is utilized for training the parameters of our neural network structure. 
Additionally, we propose a novel, adaptive, training procedure to handle state-constrained cases. 
We applied our method in path planning problems with obstacle and collision avoidance, including a system with Newtonian dynamics and non-convex environment, up to dimension 512.

The present paper extends the idea of encoding the knowledge about optimal control problems and Hamilton ODE systems in the neural network, originally presented in~\cite{meng2022sympocnet}.  
Here, to handle more general dynamics and control problems, 
we introduce 
a time-dependent module in the coordinate transformation and parameterize the latent space by the solution to a particular linear quadratic regulator (LQR) problem. Moreover, we keep the dimension of the physical space and the latent space be the same.  

This paper is organized as follows: Section~\ref{sec:bkgd} provides an overview of the background information to be utilized, covering Hamiltonian systems and symplectic maps (Section~\ref{sec:bkgd_HODE_symplectic}), SympNets (Section~\ref{sec:sympnet}), optimal control problem setup (Section~\ref{sec:problem}), SympOCNets (Section~\ref{sec:bkgd_sympocnet}). In Section~\ref{sec:optctrl_constrained}, we introduce our problem setup (Section~\ref{sec:problem_setup_1}) and the TSympOCNet architecture (Section~\ref{sec:tsympocnet}), which consists of a latent representation and a coordinate transformation. The comprehensive training algorithm is elaborated in Section~\ref{sec:training}. In Section \ref{sec:numerical_results}, we apply our method to multi-agent path planning problems with obstacle and collision avoidance. Simulations with single and four agents demonstrate robustness and effectiveness, validated against the shooting method in Section~\ref{sec:simpler}.
Sections~\ref{sec:highd} and~\ref{sec:nonconvex} showcase the effectiveness and efficiency of TSympOCNet on handling high-dimensional problems and problems with complicated constraints. Finally, a summary of findings and concluding remarks are presented in Section~\ref{sec:conclusion}.

\section{Preliminary background}\label{sec:bkgd} 

This section provides some background materials used in the remainder of the paper. 
In~\Cref{sec:bkgd_HODE_symplectic}, we provide a brief summary of symplectic maps, Hamiltonian ODEs and their relations. 
In~\Cref{sec:sympnet}, we review the SympNet architecture. In~\Cref{sec:problem}, we briefly recall the optimal control problem, the corresponding Hamiltonian system, and shooting method.  
~\Cref{sec:bkgd_sympocnet} is a review of the SympOCNet architecture, which is the original work by the authors to employ neural network architecture for solving optimal control problems. 

\subsection{Hamiltonian systems and symplectic maps} \label{sec:bkgd_HODE_symplectic}

\begin{definition}[Symplectic maps]\label{def:symp_map}
Let $U$ be an open set in $\R^{2n}$. A differentiable map $\phi:U\rightarrow\R^{2n}$ is called symplectic if the Jacobian matrix $\nabla \phi$ satisfies
\begin{equation}\label{eq:symp_map}
    \nabla \phi^T(\bz) J \nabla \phi(\bz) = J, \ \forall \ \bz\in U \ , 
\end{equation}
where $J$ is a matrix with $2n$ rows and $2n$ columns defined by
\begin{equation}\label{eqt:bkgd_def_J}
    J:=\begin{pmatrix} \bm{0} & I_{n} \\ -I_{n} & \bm{0} \end{pmatrix} \ ,
\end{equation}
and $I_n$ denotes the identity matrix with $n$ rows and $n$ columns.
\end{definition}
The Hamiltonian ODE system is a dynamical system taking the form 
\begin{equation}\label{Hh}
    \dot{\bz}(s) = J\nabla H(\bz(s)) \ ,
\end{equation}
where $J$ is the matrix defined in~\eqref{eqt:bkgd_def_J} and $H:U\rightarrow \R$ is a function called Hamiltonian.
The Hamiltonian systems and symplectic maps are highly related to each other. To be specific, the Hamiltonian structure is preserved under the change of variable using any symplectic map. This result is stated in the following theorem. 
\begin{theorem}\cite[Theorem~2.8 on p.~187]{hairer2006geometric}\label{thm:bkgd_symplectic_HODE}
Let $U$ and $V$ be two open sets in $\R^{2n}$. Let $\phi: U \rightarrow V$ be a change of coordinates such that $\phi$ and $\phi^{-1}$ are continuously differentiable functions. If $\phi$ is symplectic, the Hamiltonian ODE system $\dot{\bz}(s) = J\nabla H(\bz(s))$ can be written in the new variable $\bw = \phi(\bz)$ as
\begin{equation}\label{eq:change_var}
    \dot{\bw}(s) = J\nabla \tilde{H}(\bw(s)) \ ,
\end{equation}
where the new Hamiltonian $\tilde{H}$ is defined by 
\begin{equation}\label{eqt:thm1_def_tildeH}
    \tilde{H}(\bw) = H(\bz) = H(\phi^{-1}(\bw)),  \ \forall \ \bw\in V \ .
\end{equation}
Conversely, if \textcolor{black}{$\phi: U \rightarrow V$ is a change of coordinates that transforms every Hamiltonian system to another Hamiltonian system by \eqref{eq:change_var} and~\eqref{eqt:thm1_def_tildeH}, then $\phi$ is symplectic}.
\end{theorem}
\Cref{thm:bkgd_symplectic_HODE} indicates that a symplectic map can transform a Hamiltonian ODE system to another one, which is potentially in a simpler form. 
This is the starting point of our proposed method. 

\subsection{SympNet}\label{sec:sympnet}
SympNet is a neural network architecture proposed in~\cite{jin2020sympnets} to approximate symplectic transformations. 
There are different kinds of SympNet architecture. In this paper, we use the G-SympNet. 
For other architectures, we refer the readers to~\cite{jin2020sympnets}.

\begin{definition}
    Let $\hat{\sigma}_{K^i,\ba^i,\bb^i}(\bx):=(K^i)^\top(\ba^i\odot\sigma(K^i\bx+\bb^i))$ for any $\bx\in\R^n$, where $\sigma$ is an activation function (e.g., sigmoid, ReLU, $\cdots$), and $\odot$ denotes the componentwise multiplication. Any G-SympNet $\varphi^{(G)}$ is an alternating composition of the following two parameterized functions:
 \begin{equation}\label{eq:symp_net}
    \begin{split}
    &\mathcal{G}_{up}^{i}\begin{pmatrix} \bx \\ \bp \end{pmatrix}=\begin{pmatrix} \bx\\\bp+\hat{\sigma}_{K^i,\ba^i,\bb^i}(\bx) \end{pmatrix} \quad \forall \ \bx,\bp\in\R^n \ ,\\ &\mathcal{G}_{low}^{i}\begin{pmatrix} \bx \\ \bp \end{pmatrix}=\begin{pmatrix}  \hat{\sigma}_{K^i,\ba^i,\bb^i}(\bp)+\bx \\ \bp \end{pmatrix}\quad \forall \ \bx,\bp\in\R^n \ ,    \\
    &\varphi^{(G)} = \mathcal{G}_{up}^{N}\circ \mathcal{G}_{low}^{N}\cdots \mathcal{G}_{up}^1\circ \mathcal{G}_{low}^1 \quad \text{or} \quad \varphi^{(G)}= \mathcal{G}_{low}^N\circ \mathcal{G}_{up}^N\cdots \mathcal{G}_{low}^1\circ \mathcal{G}_{up}^1 \ ,
    \end{split}
    \end{equation}
     where the learnable parameters are the matrices $K^i\in\mathbb{R}^{l\times n}$ and the vectors $\ba^i,\bb^i\in\R^l$, for all $i \in \{1, \cdots, N\}$.
    The dimension $l$ (which is the dimension of $\ba^i,\bb^i$ as well as the number of rows in $K^i$) is a positive integer that can be tuned, called the width of SympNet. $N$ is the number of layers of SympNet.
\end{definition}
In \cite{jin2020sympnets}, it is proven that G-SympNets are universal approximators within the family of symplectic maps. Note that it is easy to obtain the inverse map of a G-SympNet, since we have explicit formulas for $(\mathcal{G}_{up}^i)^{-1}$ and $(\mathcal{G}_{low}^i)^{-1}$ given as follows
\begin{equation}\label{eq:symp_inverse}
    (\mathcal{G}_{up}^i)^{-1}\begin{pmatrix} \bx\\\bp \end{pmatrix}=\begin{pmatrix}
    \bx \\ \bp-\hat{\sigma}_{K^i,\ba^i,\bb^i}(\bx)  \end{pmatrix},\quad (\mathcal{G}_{low}^i)^{-1}\begin{pmatrix}  \bx\\\bp \end{pmatrix}=\begin{pmatrix}  \bx-\hat{\sigma}_{K^i,\ba^i,\bb^i}(\bp) \\ \bp \end{pmatrix}.
    \end{equation}

\subsection{Optimal Control Problem, Hamiltonian ODE system and shooting method}\label{sec:problem}
Let us consider a finite horizon deterministic optimal control problem 
\begin{subequations} \label{eqt:optctrl_constrained}
\begin{equation}\label{cost}
\begin{split}
 &\inf_{\bu(\cdot) \in \U }
 \left\{\int_0^T L(\bx(s),\bu(s)) ds \right\} \ , %
\end{split}
\end{equation}\label{constraint}
over the set of trajectories $(\bx(\cdot),\bu(\cdot))$ satisfying 
\begin{equation}
\left\{ \begin{aligned}
  & \dot{\bx}(s)=f(\bx(s), \bu(s)), \ \forall \  s\in [0,T] \ , \\
  &\bx(0)=\bx^{(0)}, \  \bx(T)=\bx^{(T)}  \ .
\end{aligned} \right.
\end{equation}
\end{subequations}
Here, $\U := \{\bu:[0,T] \to U \subset \R^m \mid \bu(\cdot) \text{ is measurable} \}$ is the set of controls.
The Lagrangian (or running cost) $L: \R^{n} \times U \to\R$, the dynamics $f: \R^n \times U \to \R^n$ are given functions, with basic regularity properties: continuous w.r.t all variables, and Lipschitz continuous w.r.t $x$, for every $u \in U$.  

A necessary optimality condition for the problem~\eqref{eqt:optctrl_constrained} is given by the Pontryagin's maximum principle (see for instance~\cite{trelat2005controle}). By introducing the so-called costate $\bp:[0,T] \to \R^n$,  
the optimal trajectory $\bx$ of the control problem together with the costate $\bp$ satisfy the following Hamiltonian ODE system
\begin{subequations}\label{eqt:Hamiltonian_ODE}
\begin{equation}\label{Hh_control} 
\begin{dcases}
\dot{\bx}(s) = \nabla_{\bp} H(\bx(s), \bp(s)) \ ,\\
\dot{\bp}(s) = -\nabla_{\bx} H(\bx(s), \bp(s)) \ ,\\
\end{dcases}
\end{equation}
for every $s \in [0,T]$, with initial, final condition 
\begin{equation}\label{eqt:Hamiltonian_ODE_initial}
\bx(0) = \bx^{(0)}, \ \bx(T) = \bx^{(T)} \ ,
\end{equation}
where \textcolor{black}{the Hamiltonian $H\colon \R^{2n} \to \R$ is defined by}
\begin{equation}\label{eqt:defH_unconstrained}
    H(\bx,\bp) = \max_{\bu\in U } \{ \langle \bp, f(\bx, \bu)\rangle - L(\bx,\bu)\}  \ .
\end{equation}
\end{subequations} 
Notice that~\eqref{Hh_control} is consistent with~\eqref{Hh} by taking $\bz = (\bx,\bp)$. 

One classical numerical method to solve system~\eqref{eqt:Hamiltonian_ODE} is the shooting method, which converts the initial-terminal system to a root finding problem (see for instance~\cite{pesch1994practical}). Let us consider the same Hamiltonian system~\eqref{eqt:Hamiltonian_ODE} but with~\eqref{eqt:Hamiltonian_ODE_initial} replaced by the initial condition
\begin{equation}\label{initial_only}
\bx(0) = \bx^{(0)}, \ \bp(0) = \bp^{(0)} \ .
\end{equation}
For a given $\bp^{(0)}$, solving this system gives a trajectory together with the co-state $(\bx(s),\bp(s))$, for every $s\in[0,T]$. We can then define a shooting function $S$ which maps $\bp^{(0)}$ to the value of the final conditions, for instance $S(p^{(0)}) = \|x(T)-x^{(T)} \|$. Finding the zero of $S$ gives a solution to~\eqref{eqt:Hamiltonian_ODE}. This is typically done by (quasi-)Newton method. However, this method is significant sensitive to the initial guess. The algorithm may either not converge at all, or converge to a local minimum in the presence of nonconvexity. 

\subsection{SympOCNets}\label{sec:bkgd_sympocnet}
\begin{figure}[h]
    \centering
    \includegraphics[width=0.9\textwidth]{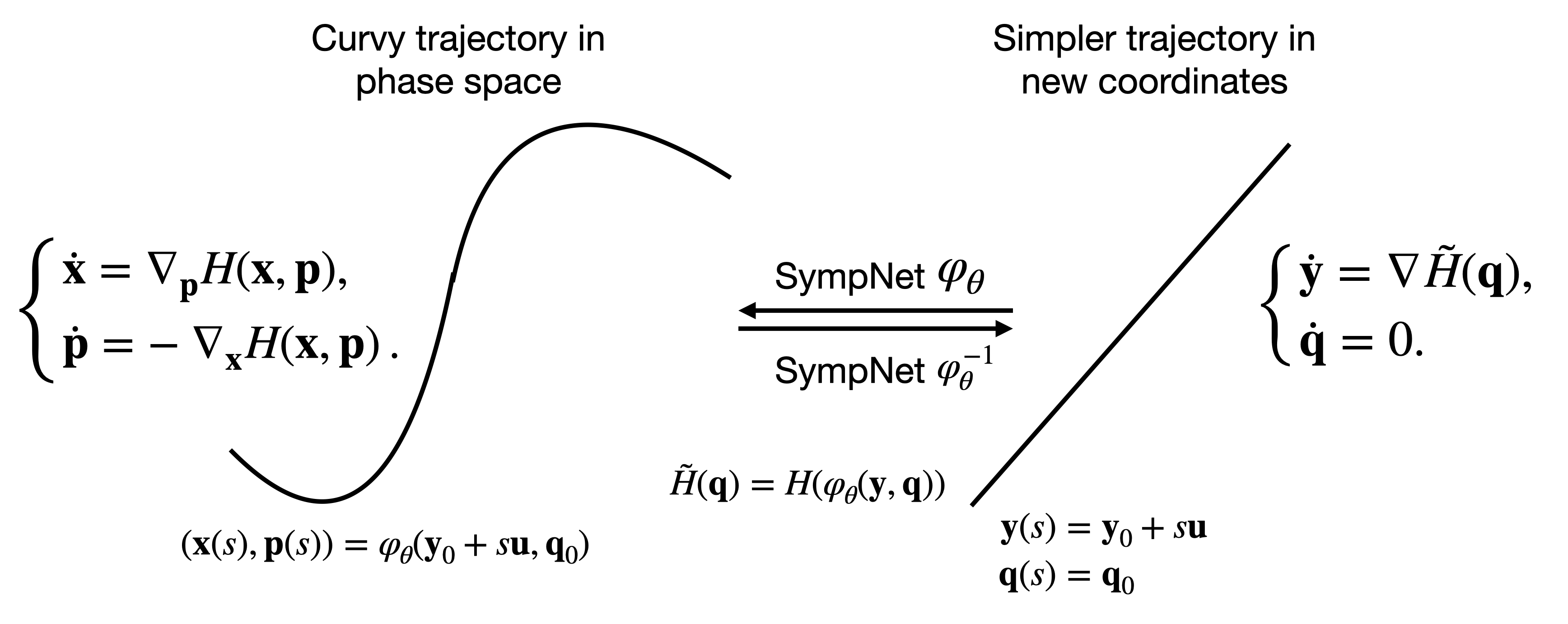}
    \caption{\textbf{An illustration of the SympOCNet method.} SympOCNet $\varphi_{\theta}$ maps curvy trajectory in the phase space to affine maps in new coordinates.}
    \label{fig:SympOCNet_orig}
\end{figure}
SympOCNet was introduced in \cite{meng2022sympocnet} to first solve the optimal control problem~\eqref{eqt:optctrl_constrained} together with the system~\eqref{eqt:Hamiltonian_ODE}, and then extended to include a state constrained case. It focuses on a special case of dynamics
\begin{equation}\label{special_dynamics}
\dot{\bx}(s) = \bu(s), \ \forall \ s \in [0,T] \ .
\end{equation}
The crucial step of SympOCNet is employing a change of variables technique in the phase space, by a symplectic map represented by SympNet architecture, to transform the original Hamiltonian ODE system~\eqref{eqt:Hamiltonian_ODE} into a simper one that can be easily solved.
In more detail, assuming given a symplectic map $\phi :\R^{2n} \to \R^{2n}$, the change of variables is defined by $(\by ,\bq) = \phi(\bx,\bp)$. 
Then, if $(\bx(s), \bp(s))$, for $s\in[0,T]$, is a solution of system~\eqref{eqt:Hamiltonian_ODE}, we assume $(\by(s), \bq(s))$, for $s\in[0,T]$, is a solution of a new Hamiltonian ODE system, with the Hamiltonian taking the form
\begin{equation}\label{SympOCNet_H}
\Tilde{H}(y,q) = H(\phi^{-1}(y,q)) \ .
\end{equation}
For computational efficiency, it is also assumed in~\cite{meng2022sympocnet} that the new Hamiltonian $\tilde{H}$ does not depend on $\by$, leading to a system:
\begin{subequations}\label{eqt:Hamiltonian_ODE_so}
\begin{equation} 
\begin{dcases}
\dot{\by}(s) = \nabla_{\bp} \Tilde{H}( \bq(s)) \ ,\\
\dot{\bq}(s) = 0  \ ,\\
\end{dcases}
\end{equation}
for every $s\in[0,T]$, with initial and terminal condition
\begin{equation}
(\phi^{-1})^{(1)}(\by(0),\bq(0)) = \bx^{(0)}, \ (\phi^{-1})^{(1)}(\by(T),\bq(T)) = \bx^{(T)} \ ,
\end{equation}
\end{subequations}
where $f^{(1)}$ denotes the first $n$ output components of any $f:\R^{2n} \to \R^{2n}$. 
The inverse of this symplectic map $\phi^{-1}$ is then approximated through a parameterized family of symplectic maps $\varphi_{\theta}$, where $\theta$ represents the unknown parameters to be learned through SympNet architecture. 
The solution to the original problem can be obtained by mapping the trajectory back to original phase space through $\varphi_{\theta} = \phi^{-1}$.
We refer to~\cite{meng2022sympocnet} for more details.

\section{TSympOCNet for optimal control problems with state constraints}\label{sec:optctrl_constrained}
In this section, we purpose to solve optimal control problems with time-dependent symplectic optimal control network (TSympOCNet), a neural network parameterized by an LQR latent representation and a time-dependent symplectic coordinate transformation.
\subsection{Control Problem and Sketch of the Architecture}\label{sec:problem_setup_1}
\subsubsection{Problem Setup}\label{sec:problem_setup}
We are interested in the optimal control problem~\eqref{eqt:optctrl_constrained} with more general, potentially nonlinear with respect to $x$, dynamics
\begin{equation}\label{eq:oc_nonlinear}
\dot{\bx}(s) = f(\bx(s)) + B\bu(s), \ \forall \ s \in [0,T] \ .
\end{equation}
Moreover, we consider the state constraint case, that is in~\eqref{constraint}, we further require
\begin{equation}\label{state_constraint}
h(\bx(s)) \geq 0, \ \forall \ s \in [0,T] \ .
\end{equation}
Here, $h: \R^n \to \R^{n'}$ 
is a function imposing the constraints on the state variable $\bx$. For instance, to avoid obstacle, $h$ can be defined using the signed distance function to the obstacles. We adopt a similar technique, the soft penalty method, as in~\cite{meng2022sympocnet} to convert the constrained problem into an unconstrained one. Consider a penalty function $U_{\epsilon,l}:\R^{n'} \to \R $, dependent on positive hyperparameters $\epsilon$ and $l$, and is defined as follows: for every  $\h=(\h_1,\dots,\h_{n'})\in\R^{n'}$
\begin{subequations}
    \begin{equation}\label{eqt:def_beta_hd}
        U_{\epsilon, l}(\h) :=  \max_{i\in \{1,\cdots,{n'}\}} U_{\epsilon, l}^i(\h_i) \ ,
    \end{equation}
    where for every $i \in \{1,\dots,m \}$, $U_{\epsilon, l}^i: \R \to \R$ is defined by
    \begin{equation}\label{eqt:def_U_1d}
    U_{\epsilon, l}^i(\h_i) :=  \begin{dcases}
       -\epsilon\log(\h_i), & \text{if}\ \h_i>l \ , \\
       -\epsilon\log(\h_i) + \frac{\epsilon}{2}\left(\left(\frac{\h_i-2l}{l}\right)^2 - 1\right), & \text{if}\ \h_i\leq l \ .
   \end{dcases}
    \end{equation}
\end{subequations}
With this penalty function, the state constrained problem is converted to an unconstrained one with Lagrangian (or running cost)
\begin{equation} \label{eqt:optctrl_penalty_general}
\end{equation}
\begin{remark}\label{rm_penalty}
The induction of~\eqref{eqt:optctrl_penalty_general} indeed replaces the hard constraint $h$ by a penalty term in the cost functional of the original problem. Moreover, by an elementary computation, one can find that given two monotone sequences of parameters $\{\epsilon_l\}, \{l_k \}$ such that 
\begin{equation}\lim_{k\to \infty} \epsilon_k, l_k \to 0,\ \text{ and } \lim_{k\to \infty} \epsilon_k \log(l_k), \frac{l_k^2}{\epsilon_k} \to 0 \ ,
\end{equation}
the penalty function~\eqref{eqt:optctrl_penalty_general} tends to the indicator function of $[0,+\infty)$. Hence the new unconstrained problem is equivalent to the original problem with constraint $h(x) \geq 0$. This property in particular motivates us to develop a new adaptive training procedure, which will be detailed later.
\end{remark}

Although the method we proposed in the present paper is applicable in a wider range of problems, 
we are particularly interested in the context of path planning problems involving obstacle and collision avoidance. 
We assume that the control system consists of $M$ elementary dynamical subsystems in interaction. 
The state of each subsystem has a dimension $d_x$, 
where the concatenated state variable $\bx = (\bx_1, \cdots, \bx_M)\in\R^{Md_x}$. 
Similarly, the control vector $\bu = (\bu_1, \cdots, \bu_M)\in\R^{Md_u}$. 
For each subsystem $i$, denote by $F_i: \R^{d_x} \to \R$ an individual ``potential energy'' function associated with state $x_i$, 
and $G_i: \R^{d_u} \to \R$ a convex individual ``kinetic energy'' function associated with control $u_i$. 
The individual dynamics of the subsystem have the form $\dot{x}_i(s) = f_i(x_i(s)) + B_i u_i(s)$, 
for every $s\in[0,T]$. 
Moreover, $h$ will represent both the interaction (i.e., collision avoidance) and the state constraint (i.e., obstacle) of the problem. 
We look for a trajectory $x(s)$, for $s \in [0,T]$, minimizing the total action functional. 
This can be interpreted using the framework of optimal control problem~\eqref{eqt:optctrl_constrained}  
by taking   
\begin{equation}\label{model_oc}
\begin{split}
     &L(\bx, \bu) = \sum_{i=1}^M \Big(F_i(\bx_i) + G_i(\bu_i) \Big) \ , \\
     &\dot{\bx}_i(s) = f_i(\bx_i(s)) + B_i\bu_i(s), \  \forall \ s\in [0,T] \text{ and } \ i \in\{ 1,\cdots, M\} \ ,
\end{split}
\end{equation}
together with state constraint~\eqref{state_constraint}. 
Observe that the Lagrangian (or running cost) appearing in~\eqref{model_oc} is the sum of kinetic and potential energy, 
instead of their difference as in classical mechanics. 
Lagrangians of the form~\eqref{model_oc}, 
in which the potential is typically coercive (tending to $\infty$ as $\| x\| \to \infty$), 
are the most natural ones in the context of optimal control. 
In particular, thanks to coercivity of the potential, the minimization problem is well defined over an arbitrary time horizon. 

Combing the techniques of converting the state constraints to running cost as in~\eqref{eqt:optctrl_penalty_general}, we have that the Hamiltonian of the problem takes the form   
\begin{equation}\label{eq:hamiltonian2}
\begin{split}
    H_{\epsilon, l}(\bx,\bp) = &\max_{\bu } \left\{ \langle\bp,f(\bx)+B\bu\rangle-L(\bx, \bu) - U_{\epsilon, l}(h(\bx)) \right\}\\
    = &\sum_{i=1}^M\Big(\langle \bp_i,f_i(\bx_i)\rangle - F_i(\bx_i)+ G_i^*(B^T_i\bp_i)\Big) - U_{\epsilon, l}(h(\bx)) \ ,
\end{split}
\end{equation}
where $G_i^*$ denotes the Legendre-Fenchel transform of $G_i$. 
\subsubsection{Motivation for the New Architecture}
In SympOCNet architecture, for computational efficiency, we assume that $\tilde{H}$ does not depend on $\by$ in~\eqref{eqt:Hamiltonian_ODE_so}, leading to simpler affine maps trajectories in the new coordinates (see~\Cref{fig:SympOCNet_orig}). However, the existence of such a symplectic map that transforms~\eqref{eqt:Hamiltonian_ODE} to~\eqref{eqt:Hamiltonian_ODE_so} is not guaranteed, in particular since $\phi$ is itself an isomorphism. 

In the present paper, we aim to handle more general dynamics of the form~\eqref{eq:oc_nonlinear}, then it is natural to consider a more expressive Hamiltonian system in the latent space. So we consider a Hamiltonian ODE system in the latent space, with the Hamiltonian $\Tilde{H}$ depending both on the state and co-state. Moreover, for a Hamiltonian system of the form~\eqref{Hh}, we denote by $\Phi_t$ the map sending any initial condition $z_0$ to $z(t)$, for a fixed $t$, that is the solution map. Recall that if $H$ is of class $\Cc^2$, $\Phi_t$ is a symplectic transformation (see for instance~\cite{hairer2006geometric}).   

\begin{figure}[h]
    \centering
    \includegraphics[width=0.9\textwidth]{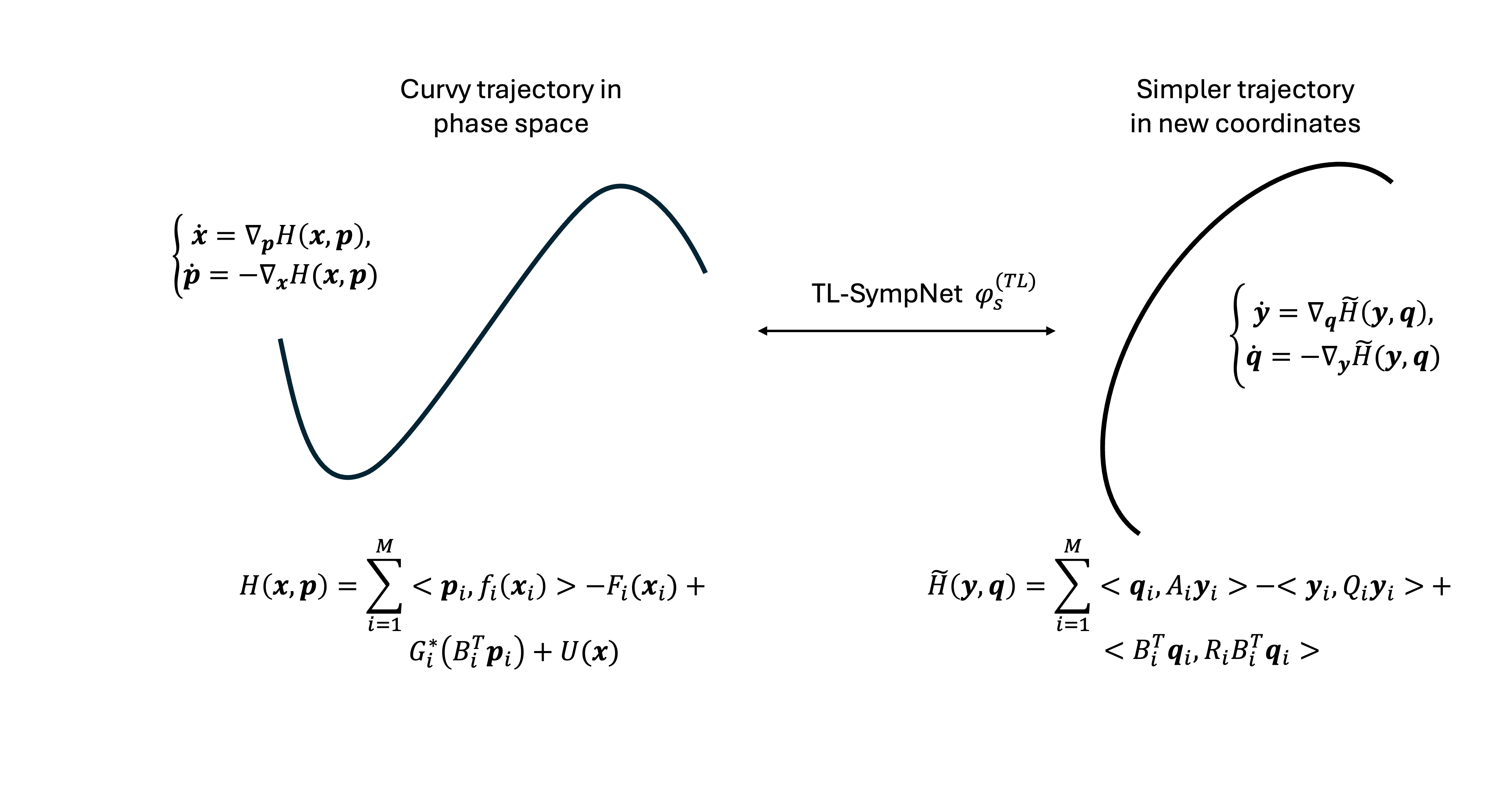}
    \caption{\textbf{An illustration of the time-dependent SympOCNet method.} We propose to use a time-dependent SympNet $\varphi_s$ to map a curvy trajectory in the phase space to a simpler trajectory in new coordinates and solve the corresponding Hamiltonian system of the optimal control problem.}
    \label{fig:SympOCNet}
\end{figure}

Based on these observations, it is tempting to assume that, by parameterizing the coordinate transformation through a time-dependent symplectic map, we can transform the original Hamiltonian ODE system into a simper one. 
Then, the solution of the original system at a fixed time $t$ is the pre-image of the solution in the latent system by this time-dependent symplectic map. 
A sketch of the present structure is shown in~\Cref{fig:SympOCNet}, and 
it will be detailed in the following sections.

\subsection{TSympOCNet architecture}\label{sec:tsympocnet}
\subsubsection{Latent representation}
In the latent space, we consider a simpler Hamiltonian ODE system, with Hamiltonian 
\begin{subequations}\label{system_latent}
\begin{equation}\label{H_latent}
    \Tilde{H}(\by, \bq) = \sum_{i=1}^M\Big(\langle \bq_i,A_i\by_i\rangle - \langle\by_i, Q_i\by_i\rangle + \langle B^T_i\bq_i, R_i B^T_i\bq_i\rangle\Big) \ ,
\end{equation}
where $A_i = \textbf{J}f_i(0)$, $Q_i = \frac{1}{2}\textbf{H}(F_i)(0)$, $R_i = \frac{1}{2}\textbf{H}(G_i^*)(0)$. Here $\textbf{J}f_i(0)$ represents the Jacobian of $f_i$ at 0 and $\textbf{H}(F_i)(0)$ represents the Hessian of $F_i$ at 0. 
The Hamiltonian system then has the form
\begin{equation}\label{Hsystem_latent}
\begin{dcases}
\dot{\by_i}(s) = \nabla_{\bq_i} \Tilde{H}(\by(s), \bq(s)), \ &  i  \in \{1,\dots, M\} \ ,\\
\dot{\bq_i}(s) = -\nabla_{\by_i} \Tilde{H}(\by(s), \bq(s)) , \ & i \in \{1,\dots, M \} \ ,  \\
\end{dcases}
\end{equation}
for every $s\in [0,T]$. Moreover, we set the same initial, final condition as the original Hamiltonian ODE system, i.e., 
\begin{equation}\label{if_latent}
y(0) = x^{(0)}, \ y(T) = x^{(T)} \ .
\end{equation}
\end{subequations}
Notice that there is a one-to-one correspondence between the latent Hamiltonian system~\eqref{system_latent} and the following linear quadratic regulator (LQR) problem
\begin{equation} \label{LQR}
\begin{aligned}
   &\inf_{\boldsymbol{u}(\cdot)} \int_0^T \sum_{i=1}^M \left(  \boldsymbol{y}_i(s)^T Q_i \boldsymbol{y}_i(s) + \boldsymbol{u}_i(s)^T R_i \boldsymbol{u}_i(s) \right) ds \\
    &\text{s.t.} 
    \begin{cases} 
         \dot{\boldsymbol{y}_i}(s) = A_i \boldsymbol{y}_i(s) + B_i \boldsymbol{u}_i(s), \quad \forall i \in\{ 1,\cdots, M\} \text{ and }  \forall s\in [0,T] \ , \\
         \boldsymbol{y}(0) = \boldsymbol{x}^{(0)}, \  \boldsymbol{y}(T) = \boldsymbol{x}^{(T)} \ . \\
    \end{cases}
\end{aligned}
\end{equation}

The Hamiltonian system~\eqref{system_latent} can be formulated as
\begin{equation} \label{eqt: LQR Hamiltonian ode}
\begin{aligned}
    \begin{cases} 
         \displaystyle \frac{d}{ds} \begin{pmatrix} \by_i(s) \\ \bq_i (s) \end{pmatrix} = \begin{pmatrix}
             A_i & - B_i R_i^{-1} B_i^T \\ -Q_i & -A_i^T
         \end{pmatrix} \begin{pmatrix} \by_i(s) \\ \bq_i (s) \end{pmatrix}, \ i \in\{ 1,\cdots, M\} \ ,\\
    \end{cases}
\end{aligned}   
\end{equation}
for every $s\in[0,T]$, with initial, final condition~\eqref{if_latent}. 
We call the dynamic matrix in \eqref{eqt: LQR Hamiltonian ode} the Hamiltonian matrix, and denote as $H_i$. 
The solution of this two-point boundary value system can be explicitly obtained. 
Notice that $(\by(s),\bq(s)) = ((\by_i)(s),\bq_i(s))_{i = \{1,\dots, M\}}$, with
\begin{equation}\label{eq:linear_hamilton}
    \begin{pmatrix} \by_i(s) \\ \bq_i (s) \end{pmatrix} = e^{H_is} \begin{pmatrix} \bx^{(0)}_{i} \\ \bq^{(0)}_{i} \end{pmatrix} \ ,
\end{equation}
solves the system~\eqref{eqt: LQR Hamiltonian ode}  with initial condition $\by(0) = \bx^{(0)}, \bq (0)  = \bq^{(0)}$. So if we can find an appropriate $\bq_0$ such that the solution also satisfies $\by(T)=\bx^{(T)}$, it is also a solution of the boundary value problem \eqref{eqt: LQR Hamiltonian ode}. 
Indeed, solution of the linear equations 
\begin{equation}\label{eq:init_q0}
    \by_i(T)  = \begin{pmatrix} I & 0 \end{pmatrix} \cdot e^{H_iT} \begin{pmatrix} \bx_{0i} \\ \bq_{0i} \end{pmatrix} 
\end{equation}
provides the desirable $\bq_0$.

\begin{remark} 
The method proposed above to solve the Hamiltonian ODE system in the latent space has a computational complexity $O(Md_x)$, in the sense of algorithmic operation. 
It is worth noticing that, on the other hand,  
the problem \eqref{LQR} can also be approximated through a finite horizon LQR problem with quadratic terminal cost
\begin{equation} \label{LQR terminal cost}
\begin{aligned}
   &\inf_{\bu(\cdot)} \int_0^T   \by(s)^T Q \by(s) + \bu(s)^T R \bu(s) ds + \by(T)^T Q_f \by(T)\\
    &\text{s.t.} 
    \begin{cases} 
         \dot{\by}(s) = A \by(s) + B \bu(s), \quad \forall s\in [0,T] \ , \\
         \by(0) = \bx \ .
    \end{cases}
\end{aligned}
\end{equation}
Let us define the value function $V$ which associates with any $(x,t) \in \R^{Md_x} \times [0,T]$ the infimum of $\int_t^T \by(s)^T Q \by(s) + \bu(s)^T R \bu(s) ds + \by(T)^T Q_f \by(T) $, under the constraints in~\eqref{LQR terminal cost}. Then $V(x,t)$ will remain quadratic form $x^T P(t) x$ for every $t \in [0,T]$, with $P(t)$ the solution of the following differential Riccati equation 
\begin{equation} \label{eq: Riccati_P}
\left\{
\begin{aligned}
    \dot{P}(s) &= -A^T P(s) - P(s) A + P(s) B R^{-1} B^T P(s) - Q,
    \ \forall \ s \in [t,T] \ , \\
    P(T) &= Q_f \ .
\end{aligned}
\right.
\end{equation}
Notice that solving the Riccati equation~\eqref{eq: Riccati_P}, although it has a computational complexity $O(Md_x^3)$, gives an optimal control as a function of the state $\bx$. 
It is also called feedback control or closed loop control, which is known to have advantages in real application, for instance the solution is robust against system perturbations. 
Constructing a robust neural network architecture for the control problem is of independent interest, and we leave it as a direction for future work. 
\end{remark}

%

\subsubsection{Coordinate transformation}\label{sec:tl_transform}
In the reminder of this paper, we will denote by $(\bx(s), \bp(s))_{s\in[0,T]}$ the solution of the Hamiltonian system in the phase space with Hamiltonian~\eqref{eq:hamiltonian2}, and $(\by(s),\bq(s))_{s \in [0,T]}$ the solution of the Hamiltonian system~\eqref{system_latent} in the latent space. The goal is to represent and approximate the inverse of an unknown time-dependent symplectic map $\varphi$, that can transform $(\bx(t),\bp(t))$ to $(\by(t),\bq(t))$, for any fixed $t \in [0,T]$. 
\begin{definition}
Let $U$ be an open set in $\mathbb{R}^{2n}$. 
We call a map $\varphi\in \Cc^{1}(U\times[0,T];\R^{2n})$ time-dependent symplectic 
if for any fixed $t \in [0,T]$, $\varphi(\cdot,t)$ is symplectic. 
That is, for any $t \in [0,T]$, 
\begin{equation}
\nabla_{\bz}\varphi(\bz, t)^TJ\nabla_{\bz}\varphi(\bz, t) = J, \ \forall \  \bz\in U \ .
\end{equation}
\end{definition}
For simplicity, we use $\varphi_t$ to denote $\varphi(\cdot,t)$, for a fixed $t\in[0,T]$. 
In this paper, we consider the \textit{linear} and \textit{affine} time-dependent symplectic map $\varphi$.
\begin{remark}
We refer to a time-dependent symplectic transformation that is linear (affine) in $z$ as a linear (affine) time-dependent symplectic transformation. In the linear case, each $\varphi_t$ is a linear symplectic transformation that can be represented by a symplectic matrix. 
\end{remark}
\begin{lemma}[\cite{jin2022optimal}]\label{lemma_tsymp}
For any symplectic matrix $K\in\mathbb{R}^{2n\times2n}$, there exists $k = 5$ symmetric matrices $\{S_i\}_{i=1}^k$, such that
\begin{equation}
    K =
    \begin{pmatrix}
    I & 0 \\ 
    S_1 & I
    \end{pmatrix}
    \begin{pmatrix}
    I & S_2 \\ 
    0 & I
    \end{pmatrix} 
    \cdots
    \begin{pmatrix}
    I & 0 \\ 
    S_k & I
    \end{pmatrix}.
\end{equation}
\end{lemma}
\Cref{lemma_tsymp} shows that any linear symplectic transformation can be written as the composition of at most $k = 5$ unit block lower/upper triangular matrices. Inspired by the lemma above, we propose the TL-SympNet, which resembles the G-SympNet to represent arbitrary affine time-dependent symplectic transformation.

\begin{definition}
    For $i \in \{1,2,\dots, N \}$, let $\Tilde{\sigma}_{K^i,\ba^i,\bb^i}(\bx;t):=(K^i)^\top(\ba^i(t)\odot(K^i\bx+\bb^i))$ for any $\bx\in\R^n$, $t\in[0, T]$. Any TL-SympNet $\Tilde{\varphi}_t$ is an alternating composition of the following two parameterized functions:
    \begin{subequations}\label{eq:tl_sympnet}
 \begin{equation}\label{tlow}
    \begin{split}
    &\mathcal{T}^{low, i}_t\begin{pmatrix} \bx \\ \bp \end{pmatrix}=\begin{pmatrix} \bx\\\bp+\Tilde{\sigma}_{K^i,\ba^i,\bb^i}(\bx;t) \end{pmatrix} \quad \forall \bx,\bp\in\R^n,\\ 
    &\mathcal{T}^{up, i}_t\begin{pmatrix} \bx \\ \bp \end{pmatrix}=\begin{pmatrix}  \Tilde{\sigma}_{K^i,\ba^i,\bb^i}(\bp;t)+\bx \\ \bp \end{pmatrix}\quad \forall \bx,\bp\in\R^n,    \\
    \end{split}
    \end{equation}
    such that
    \begin{equation}
    \begin{split}
    &\Tilde{\varphi}_t = \mathcal{T}^{up,N}_t\circ \mathcal{T}^{low,N}_t\cdots \mathcal{T}^{up,1}_t\circ \mathcal{T}^{low,1}_t \text{ or } \Tilde{\varphi}_t= \mathcal{T}^{low,N}_t\circ \mathcal{T}^{up,N}_t\cdots \mathcal{T}^{low,1}_t\circ \mathcal{T}^{up,1}_t,
    \end{split}
    \end{equation}
    \end{subequations}
     where the learnable parameters are from the matrices $K^i\in\mathbb{R}^{l\times n}$, the vectors $\bb^i\in\R^l$, and fully connected neural network $\ba^i:\R\to\R^l$.
    The dimension $l$ (which is the dimension of $\bb^i$ as well as the number of rows in $K^i$) is a positive integer that can be tuned, called the width of TL-SympNet. $N$ is the number of layers of TL-SympNet. We refer to the number of layers and width of each $\ba^i$ as the sublayers and subwidth of TL-SympNet.
\end{definition}
Then, following directly from the universal approximation theorem of standard neural networks, we have the following result:
\begin{theorem}
    The map $\Tilde{\varphi}_t$ defined in~\eqref{eq:tl_sympnet} is an universal approximator within the set of affine time-dependent symplectic transformations.
\end{theorem}
Recall that the original Hamiltonian system and the Hamiltonian system~\eqref{system_latent} 
in latent space 
share the same set of boundary values $\bx^{(0)}$ and $\bx^{(T)}$. This inspires us to design a neural network architecture which preserves the boundary values of $\bx$.
\begin{definition}
    Let $\hat{\sigma}_{K^i,\ba^i,\bb^i}(\bx;t):=(K^i)^\top(t(T-t)\ba^i(t)\odot(K^i\bx+\bb^i))$ for any $\bx\in\R^n$, $t\in[0, T]$. Any boundary-value preserving TL-SympNet $\hat{\varphi}_t$ is an alternating composition of $\mathcal{T}^{low}$, defined as in~\eqref{tlow}, and a new $\hat{\mathcal{T}}^{up}$, defined as
    \begin{subequations}\label{eq:tl_sympnet2}
 \begin{equation}
 \begin{split}
      &\hat{\mathcal{T}}^{up, i}_t\begin{pmatrix} \bx \\ \bp \end{pmatrix}=\begin{pmatrix} \bx+\hat{\sigma}_{K^i,\ba^i,\bb^i}(\bp;t) \\\bp\end{pmatrix} \quad \forall \bx,\bp\in\R^n \ , \\
      \end{split}
      \end{equation}
      such that 
      \begin{equation}
      \begin{split}
    &\hat{\varphi}_t = \hat{\mathcal{T}}^{up, N}_t\circ \mathcal{T}^{low, N}_t\cdots \hat{\mathcal{T}}^{up, 1}_t\circ \mathcal{T}^{low, 1}_t \quad \text{or} \quad \hat{\varphi}_t= \mathcal{T}^{low, N}_t\circ \hat{\mathcal{T}}^{up, N}_t\cdots \mathcal{T}^{low, 1}_t\circ \hat{\mathcal{T}}^{up, 1}_t .
 \end{split}
\end{equation}
\end{subequations}
Note that $\hat{\varphi}_0^{(1)}(\bx, \bp) = \hat{\varphi}_T^{(1)}(\bx, \bp) = \bx$. Here $f^{(1)}$ denotes the first $n$ output components of any $f:\R^{2n} \to \R^{2n}$.
\end{definition}
We will use the boundary-value preserving TL-SympNet $\hat{\varphi}_t$ as the coordinate transformation module in TSympOCNet.

\subsection{The TSympOCNet Method}\label{sec:training}
In this subsection, we describe our algorithm of using TSympOCNet to solve the original optimal control problem together with the Hamiltonian system. 
We use $\hat{\varphi}^\theta_t$ instead of $\hat{\varphi}_t$  to denote the coordinate transformation (i.e., boundary-value preserving TL-SympNet), to highlight the dependence on the trainable variable $\theta = (K^i, \ba^i, \bb^i)_{i\in\{1,\dots,N\}}$. 

Let us first discretize the time horizon by $\bar{N} = \frac{T}{ \Delta t }$ steps, and denote $\bt_k = k \Delta t$, for $k \in \{0,1, \dots, \bar{N}\}$. 
The optimal trajectory in the latent space at every $\bt_k$, $(\by(\bt_k),\bq(\bt_k))$, can be computed using \eqref{eq:linear_hamilton} and \eqref{eq:init_q0}. 
The inverse of $(\varphi_t)_{t \in [0,T]}$, which maps the trajectory and Hamiltonian system of the latent space to the phase space, will be approximated by $(\hat{\varphi}^{\theta}_t)_{t\in[0,T]}$, and evaluated at 
data points $\{(\by(\bt_k), \bq(\bt_k) \}$. 

To learn the parameters $\theta$, we perform iterative steps to construct and minimize a physics-informed loss function of $\theta$. 
In each iteration, 
we randomly sample $\Tilde{N} $ points among $[0,T]$, and denote the value of them by $\{t_i\}_{i \in \{1,2,\dots, \Tilde{N} \} }$. 
Then, if $t_i \in [\bt_k, \bt_{k+1}[$, for some $k \in \{0,1,\dots, \bar{N}-1\}$, we calculate the value $(\by(t_i), \bq(t_i))$ by a linear interpolation of $(\by(\bt_k),\bq(\bt_k))$ and $(\by(\bt_{k+1}), \bq(\bt_{k+1}))$. 
Moreover, by the affine property of $\hat{\varphi}^\theta_t$ as defined in~\eqref{eq:tl_sympnet2}, the value of $\hat{\varphi}^\theta_{t_i}(\by(t_i), \bq(t_i))$ can be computed by a linear interpolation of $\hat{\varphi}^\theta_{t_i}(\by(\bt_k), \bq(\bt_k))$ and $\hat{\varphi}^\theta_{t_i}(\by(\bt_{k+1}), \bq(\bt_{k+1}))$. For simplicity, we will denote by $I^{t_i}$ such linear interpolation operator, i.e., 
\begin{equation}\label{phase_tra}
\begin{aligned}
(\bx_\theta(t_i), \bp_\theta(t_i)) &= \hat{\varphi}_{t_i}^{\theta}(\by(t_i), \bq(t_i)) \\
&:= I^{t_i}\circ \hat{\varphi}^\theta_{t_i}  ( \by(\bt_k), \bq(\bt_k) ) \ . 
\end{aligned}
\end{equation}
Then, the derivative of $(\bx_\theta(s), \bp_\theta(s))$ w.r.t $s$ at $t_i$ follows the chain rule:
\begin{equation}\label{derivative_xp}
\begin{aligned}
\begin{pmatrix}
           \dot{\bx_\theta}(t_i) \\
            \dot{\bp_\theta}(t_i)
        \end{pmatrix}
&= \textbf{Jac}\hat{\varphi}_{t_i}^\theta(\by(t_i), \bq(t_i))
\begin{pmatrix}
            \dot{\by}(t_i) \\
            \dot{\bq}(t_i)
        \end{pmatrix}
        + \frac{\partial \hat{\varphi}^\theta_s(\by(t_i), \bq(t_i))}{\partial s}\mid_{t_i} \  \\
        & = I^{t_i} \circ \textbf{Jac}\hat{\varphi}_{t_i}^\theta  ( \by(\bt_k), \bq(\bt_k) ) \begin{pmatrix}
            \dot{\by}(\bt_k) \\
            \dot{\bq}(\bt_k)
        \end{pmatrix}  
        + I^{t_i} \circ \frac{\partial \hat{\varphi}^\theta_s(\by(\bt_k), \bq(\bt_k))}{\partial s}\mid_{t_i} 
        \end{aligned} \ ,
\end{equation}
where $\textbf{Jac}\hat{\varphi}^\theta_{t_i}$ denotes the Jacobian of $\hat{\varphi}^\theta_{t_i}$ w.r.t $(\by,\bq)$. 
Notice that in implementation, $\textbf{Jac}\hat{\varphi}^\theta_{t_i}$ and $\frac{\partial \hat{\varphi}^\theta_s}{\partial s}$ are computed via automatic differentiation, and the overall computation can be further boosted via the JVP (Jacobian vector product) method in deep learning libraries, e.g., JAX. The physics-informed loss function w.r.t $\theta$ is then defined as   
 \begin{equation}\label{PI_loss}
    \begin{split}
        \mathcal{L}(\theta) = & \sum_{i=0}^{N_t-1}||\dot{\bx}_\theta(t_i) -\nabla_pH_{\epsilon, l}(\bx_\theta(t_i), \bp_\theta(t_i))|| + \\
        & \sum_{i=0}^{N_t-1}||\dot{\bp}_\theta(t_i)+\nabla_xH_{\epsilon, l}(\bx_\theta(t_i), \bp_\theta(t_i))|| \ .
    \end{split}
    \end{equation}
Then, in each iteration, we update $\theta$ by minimizing~\eqref{PI_loss}.  The optimization is performed with stochastic gradient descent based methods, e.g., Adam \cite{adam2015}. The iteration stops when the convergence is achieved, that is the loss function $\mathcal{L}(\theta)$ is less than a sufficiently small threshold $\bar{\varepsilon}$. 

To well handle the state constrain cases, motivated by~\Cref{rm_penalty}, we also perform iterative procedures to update the hyperparameters $\epsilon$ and $l$, defined in the penalty function \eqref{eqt:optctrl_penalty_general} to convert the problem to an unconstrained case. 
We start with fixed $\epsilon = \epsilon_0, l = l_0$. Then, we do the iterations of sampling in time horizon, computing~\eqref{derivative_xp}, constructing and minimizing the loss function~\eqref{PI_loss}. We then iteratively update $\epsilon,l$ via \begin{equation}
    \begin{aligned}
        &\epsilon_{j+1} = n_1 \epsilon_{j}\\
        &l_{j+1} = n_2 l_{j} \ .
    \end{aligned}
    \end{equation}
In this step, we do a fixed number, $N^{ite}$, of iterations.

The complete Algorithm is shown in~\Cref{algo}. In particular, the choice of hyperparameters $\bar{\epsilon}, n_1, n_2$ may depend on the problem, which will be detailed in the next section.

\begin{algorithm}
\caption{TSympOCNet}
\label{algo}
\begin{algorithmic}[1]
\State {Parameterize the transformation by $ \theta = (K^i, \ba^i, \bb^i)_{i\in\{1,\dots,N\}}$ as in~\eqref{eq:tl_sympnet2}.}
\State {Discretize the time horizon by $\bar{N} = \frac{T}{\Delta t}$ steps, compute $\{\by(\bt_k),\bq(\bt_k)\}_{\bt_k \in \{0,\Delta t, \dots, T\}}$ in latent space by~\eqref{eq:linear_hamilton} and \eqref{eq:init_q0}.}
\State {Set $\epsilon = \epsilon_0, l = l_0, k=1$.}
\For { $k \leq N^{ite}$}
\While {$\mathcal{L}(\theta) \geq \bar{\varepsilon}$}
\State {Sample $\Tilde{N}$ points in $[0,T]$: $0=t_0 \leq t_1 \leq \dots,t_{\Tilde{N}} = T$.}
\State {Compute the phase trajectory $(\bx_\theta(t_i), \bp_\theta(t_i))$ by~\eqref{phase_tra}.}
\State {Compute the derivative of phase trajectory w.r.t $s$ by~\eqref{derivative_xp}.}
\State {Minimize the physics-informed loss function~\eqref{PI_loss}}
\EndWhile
\State {Set $\epsilon = n_1 \epsilon$, $l = n_2 l, k= k+1$. }
\EndFor
\end{algorithmic} 
\end{algorithm}

We also present a summary of key steps of our algorithm, in the following.
\begin{enumerate}
    \item Set $s_i = i\frac{T}{\bar{N}}, i\in\{0,\cdots \bar{N}\}$. Precompute
    $\{(\by(s_i), \bq(s_i))\}_{i=0}^{N_s-1}$ using \eqref{eq:linear_hamilton} and \eqref{eq:init_q0}. 
    \item Initialize hyperparameters $\epsilon = \epsilon_0, l = l_0$ in \eqref{eqt:optctrl_penalty_general}.
    \item Uniformly sample $0 = t_0 \leq \cdots \leq t_{\Tilde{N}} = T$. Calculate $\{(\by(t_i), \bq(t_i))\}_{i=1}^{\Tilde{N}}$ by interpolation. Compute the phase trajectory $\{(\bx_\theta(t_i), \bp_\theta(t_i))\}_{i=1}^{\Tilde{N}}$ via 
    \begin{equation*}
        (\bx_\theta(t_i), \bp_\theta(t_i)) = I^{t_i}\circ \hat{\varphi}^\theta_{t_i}  ( \by(\bt_k), \bq(\bt_k) ) \ .
    \end{equation*}
    \item Compute \begin{equation*}
       \begin{aligned}
\begin{pmatrix}
           \dot{\bx_\theta}(t_i) \\
            \dot{\bp_\theta}(t_i)
        \end{pmatrix}  = I^{t_i} \circ \textbf{Jac}\hat{\varphi}_{t_i}^\theta  ( \by(\bt_k), \bq(\bt_k) ) \begin{pmatrix}
            \dot{\by}(\bt_k) \\
            \dot{\bq}(\bt_k)
        \end{pmatrix}  
        + I^{t_i} \circ \frac{\partial \hat{\varphi}^\theta_s(\by(\bt_k), \bq(\bt_k))}{\partial s}\mid_{t_i} 
        \end{aligned} \ ,
    \end{equation*}
    where $\textbf{Jac}\varphi_{t_i}^\theta$ and $\frac{d\varphi^\theta_t}{dt}$ are computed via automatic differentiation. 
    The overall computation can be further boosted via the JVP (Jacobian vector product) method in deep learning libraries, e.g., JAX.
    \item Update $\theta$ by optimizing the physics-informed loss 
    \begin{equation*}
    \begin{split}
        \mathcal{L}(\theta) = & \sum_{i=0}^{N_t-1}||\dot{\bx}_\theta(t_i) -\nabla_pH_{\epsilon, l}(\bx_\theta(t_i), \bp_\theta(t_i))|| + \\
        & \sum_{i=0}^{N_t-1}||\dot{\bp}_\theta(t_i)+\nabla_xH_{\epsilon, l}(\bx_\theta(t_i), \bp_\theta(t_i))|| \ .
    \end{split}
    \end{equation*}
    The optimization can be done with stochastic gradient descent based methods, e.g., Adam \cite{adam2015}.
    \item Repeat step 3, 4, 5 until convergence ($\mathcal{L}(\theta) < \bar{\varepsilon}$). 
    \item Update $\epsilon, l$ via
    \begin{equation*}
    \begin{split}
        \epsilon_{j+1} = n_1 \epsilon_{j}\\
        l_{j+1} = n_2 l_{j} \ .
    \end{split}
    \end{equation*}
    \item Repeat step 3, 4, 5, 6, 7 in $N^{ite}$ times.
    \item Predict the rolled-out phase trajectory $(\bx_\theta(s_i), \bp_\theta(s_i)) = \varphi_{s_i}^{\theta}(\by(s_i), \bq(s_i))$.
\end{enumerate}

\section{Applications in path planning problems with obstacle and collision avoidance}\label{sec:numerical_results}
In this section, we apply our method to path planning problems with multiple agents. 
We assume that each agent is represented by a ball in the physical space with radius $C_{d}$. We set the state and control variable to be $\bx = (\bx_1,\dots,\bx_M)\in\R^{Md_x}$, $\bu = (\bu_1,\dots,\bu_M)\in\R^{Md_u}$, where $M$ is the number of drones, and each $\bx_j = (\bw_j, \bv_j)\in\R^{d_x}$ denotes the position and velocity of the center of each drone.
We consider the Newtonian dynamics with quadratic resistance force $\dot{\bw_i} = \bv_i$, $\dot{\bv_i} = \bu_i - k\bv_i|\bv_i|$. then the dynamics can be rewritten as
\begin{equation}
\dot{\bx_i} = \begin{pmatrix}
\bv_i\\\bu_i
\end{pmatrix} = \begin{pmatrix}
O&I\\O&O
\end{pmatrix}
\bx_i - k\begin{pmatrix}
O&O\\O&I
\end{pmatrix}
\bx_i|\bx_i| +
\begin{pmatrix}
O\\I
\end{pmatrix}\bu_i.
\end{equation}

To ensure the safety and smooth operation of the agents, we employ a constraint function $h$, comprising two components: $h_1$ for obstacle avoidance and $h_2$ for collision avoidance among agents. Let $n_o$ denote the number of obstacles, labeled as $E_1,\dots, E_{n_o}$.

The function $h_1\colon \mathbb{R}^{\frac{Md_x}{2}}\to\mathbb{R}^{Mn_o}$ is defined as follows:
\begin{equation}\label{eqt:def_h1}
h_1(\bw_1, \dots, \bw_M) = (D(\bw_1), \dots, D(\bw_M)) \quad\forall \bw_1,\dots,\bw_M\in\mathbb{R}^{\frac{d_x}{2}},
\end{equation}
where $D\colon \mathbb{R}^{\frac{d_x}{2}}\to\mathbb{R}^{n_o}$ is defined as:
\begin{equation}\label{eqt:def_d}
D(\bw) = (D_1(\bw), \dots, D_{n_o}(\bw))\quad\forall \bw\in\mathbb{R}^{\frac{dx}{2}},
\end{equation}
and each $D_j$ signifies the collision status with the $j$-th obstacle $E_j$; $D_j(\bw)<0$ denotes collision of the drone at position $\bw$ with $E_j$.

For instance, $D_j$ could be equal to the signed distance function to $E_j$. The precise definition of $D_j$ relies on the shape of $E_j$, deferred to later sections for specific examples.

The function $h_2\colon \mathbb{R}^{\frac{Md_x}{2}}\to\mathbb{R}^{M(M-1)/2}$ is the collision avoidance constraint function. Each component addresses collision avoidance between a pair of drones. Considering each drone as a ball centered at a point in $\mathbb{R}^{\frac{d_x}{2}}$, collision occurs when the distance between centers is less than the sum of their radii.

Thus, the $l$-th component of $h_2$ is:
\begin{equation}\label{eqt:def_h2_distsq}
(h_2)_l(\bw_1,\dots,\bw_M) = ||\bw_i-\bw_j|| - (2C_{d}) \quad\forall \bx_1,\dots,\bx_M\in\mathbb{R}^{\frac{d_x}{2}},
\end{equation}
where $C_d$ represents the drone's radius, and $l=i+(j-1)(j-2)/2$ for any $1\leq i < j \leq M$ denotes the corresponding constraint index for the pair $(i,j)$.

Consequently, $(h_2)_l(\bw_1,\dots,\bw_M)< 0$ indicates collision between the $i$-th and $j$-th drones, thereby facilitating collision avoidance among agents.

In all of our numerical experiments, we set $T = 10$ as the terminal time. We have $d_x = 4$, $d_u = 2$, since we consider path planning problems in 2d.

\textbf{Training Warm-up.} Apart from the training procedure listed in~\Cref{sec:training}, we offer an additional step which could potentially lead to faster convergence and reducing the possibility of converging to a local but not global optimum. We refer to this step as \textit{training warm-up}. We empirically observed that the convergence of TSympOCNet depends on the initial condition $\bx(0) = \bx^{(0)} = (\bw^{(0)}, \bv^{(0)})$ and terminal condition $\bx(T) = \bx^{(T)} = (\bw^{(T)}, \bv^{(T)})$. Suppose that TSympOCNet converges faster on another pair of $(\Tilde{\bx}_0, \Tilde{\bx}_T)$. We can convert the optimal control problem in \eqref{eq:oc_nonlinear} into a sequence of $N_{opt}$ optimization problems. At step $i$, we solve the problem
\begin{equation}\label{eq:oc_change_init}
\begin{aligned}
&\inf_{\bu(\cdot) \in \U} \int_0^T \sum_{i=1}^M F_i(\bx_i(s)) + G_i(\bu_i(s)) ds \\
&\text{s.t} \left\{
\begin{aligned}
&\bx_i(s) = \begin{pmatrix}&\bw_i(s) \\ &\bv_i(s)\end{pmatrix} \ \forall \ s \in [0,T] \ ,\\
&\dot{\bx_i}(s) = f_i(\bx_i(s)) + B_i\bu_i(s), \ \forall \ s \in [0,T] \ , \\
&h(\bw(s))\geq 0,  \ \forall \ s \in [0,T] \ , \\
&\bx(0)=\frac{i}{N_{opt}}\bx^{(0)} + \frac{N_{opt} - i}{N_{opt}}\Tilde{\bx}_0, \ \bx(T)=\frac{i}{N_{opt}}\bx^{(T)} + \frac{N_{opt} - i}{N_{opt}}\Tilde{x}_T \ ,
\end{aligned}
\right.
\end{aligned}
\end{equation}
and use the obtained coordinate transformation $\varphi_s^{(i)}$ as the initialization for the next iteration. We apply this training warm-up scheduler in the second example of~\Cref{sec:simpler} and all examples in~\Cref{sec:highd}. In our experiments, to find $(\Tilde{\bx}_0, \Tilde{\bx}_T)$, we first solve a simpler optimal control problem
\begin{equation}\label{eq:oc_simpler}
\begin{aligned}
&\inf_{\bv(\cdot) \in \mathcal{V}} \int_0^T \sum_{i=1}^M F_i(\bx_i(s)) ds \\
&\text{s.t} \left\{
\begin{aligned}
&\bx_i(s) = \begin{pmatrix}&\bw_i(s) \\ &\bv_i(s)\end{pmatrix} \ \forall \ s \in [0,T] \ ,\\
&\dot{\bw_i}(s) = \bv_i(s), \ \forall \ s \in [0,T] \ , \\
&h(\bw(s))\geq 0,  \ \forall \ s \in [0,T] \ , \\
&\bw(0)=\bw^{(0)} , \ \bw(T)=\bw^{(T)} 
\end{aligned}
\right.
\end{aligned}
\end{equation}
again with TSympOCNet to get a solution $\Tilde{\bv}$. Then we set $\Tilde{\bx}_0 = (\bw^{(0)},\Tilde{\bv}(0))$ and $\Tilde{\bx}_T = (\bw^{(T)},\Tilde{\bv}(T))$.

Simulation results are presented in the following subsections. In~\Cref{sec:simpler}, we present lower-dimensional examples and benchmark the solutions with shooting methods. In~\Cref{sec:highd}, we demonstrate the capability of the existing framework to scale to 256 dimensional problems. We further compare the total running cost, constraint violation criterion and running time cost of TSympOCNet with vanilla PINN. In~\Cref{sec:nonconvex}, we consider an example with highly nonconvex obstacles. 
In~\Cref{sec:numerical_3d}, we consider a swarm path planning problem with each agent in dimension 3. 
Video animations of these examples are available online at \href{https://github.com/zzhang222/TSympOCNet}{https://github.com/zzhang222/TSympOCNet}.

\subsection{Single / four agent with circular obstacle}\label{sec:simpler}
To prevent drone collisions with obstacles, we define the constraint function $h$ as per \eqref{eqt:def_h1}. The function $D$ in \eqref{eqt:def_d} is defined by:
\begin{equation*}
D(\bw) = ||\bw-\mathring{\bw}|| - (C_o+C_d) \quad \forall \bx\in\mathbb{R}^2,
\end{equation*}
where $\mathring{\bw}$ and $C_o$ represent the center and radius of the obstacle, respectively. In this experiment, we set $C_o = 0.15$, $C_d = 0.05$.

\textbf{Single agent.} Assume $M = 1$, $L(\bx, \bu) = \frac{1}{2}\|\bu\|^2$. It is assumed that the agent starts from the top left with zero initial velocity and ends at the lower right with zero terminal velocity, i.e. $\bx^{(0)} = (-0.5, 0.5, 0, 0)$ and $\bx^{(T)} = (0.5, -0.5, 0, 0)$.

\begin{figure}[h]
    \centering
    \includegraphics[width=0.8\textwidth]{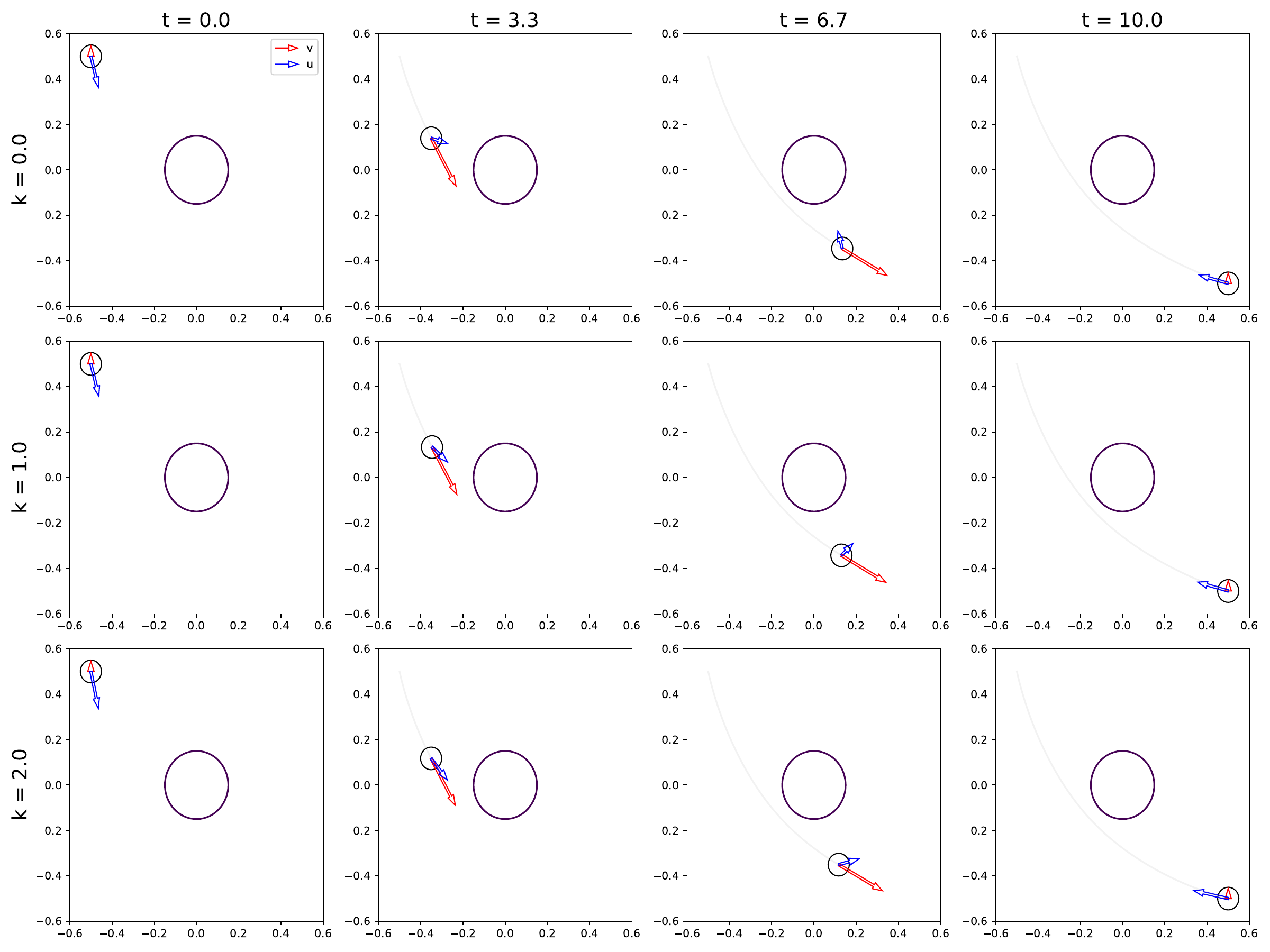}
    \caption{\textbf{The planned path and control at resistance coefficient $k = 0, 1, 2$ and time step $t = 0, 3.3, 6.7, 10$.} As resistance coefficient $k$ increases, the control signal $\bu$ exhibits prolonged alignment with the velocity $\bv$.}
    \label{fig:1d_traj}
\end{figure}

We compute the solutions to the optimal control problem \eqref{eq:oc_nonlinear} with resistance coefficient $k = 0, 1, 2$. In~\Cref{fig:1d_traj}, we plot the planned trajectory by TSympOCNet $\bx = (\bw, \bv)$ together with the control signal $\bu$. It can be seen that the same planned trajectory $\bx$ was found with varying level of $k$. As resistance coefficient $k$ increases, $\bu$ exhibits prolonged alignment with the velocity $\bv$, consistent with physical principles.

\begin{figure}[h]
    \centering
    \includegraphics[width=0.8\textwidth]{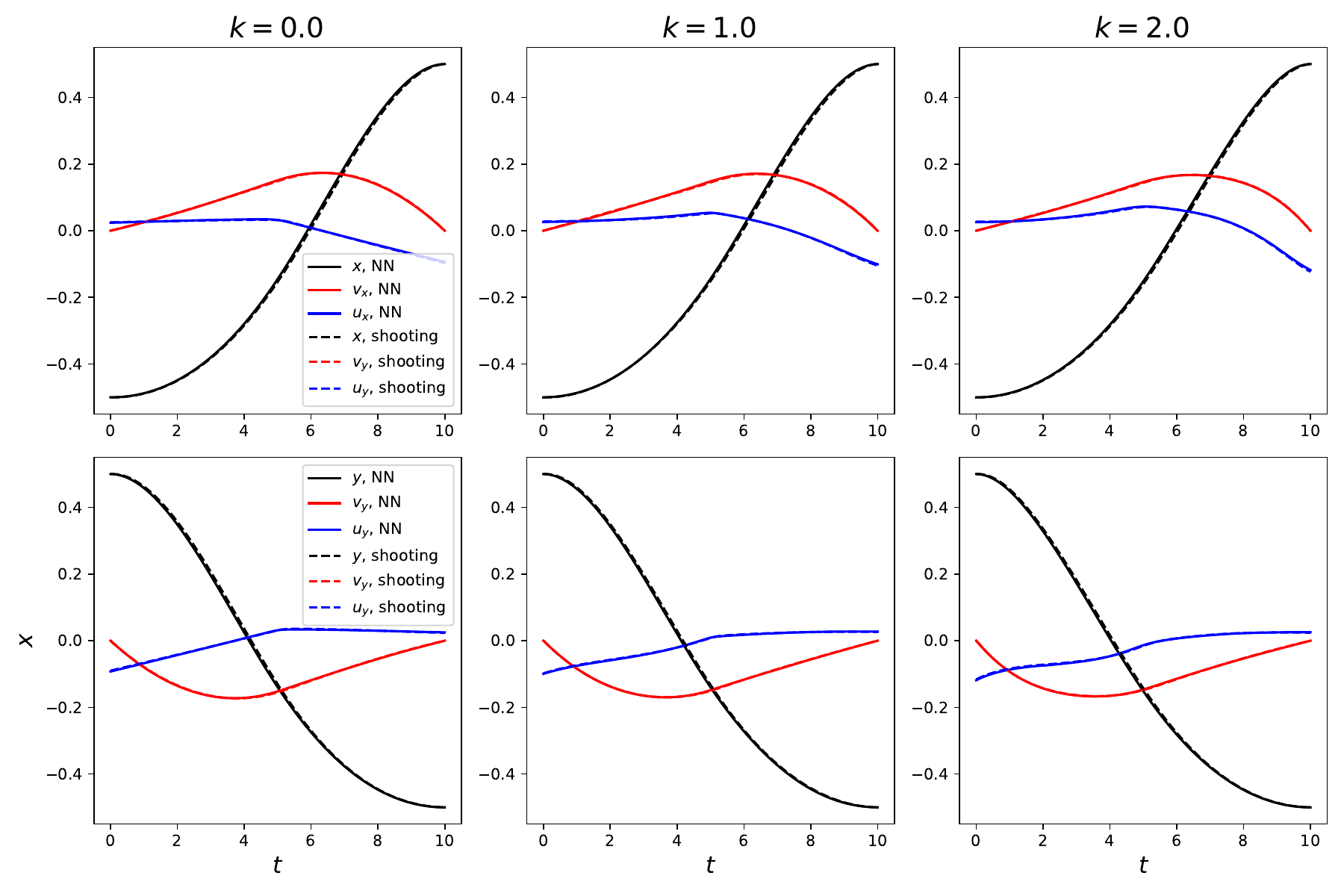}
    \caption{\textbf{The solution generated by TSympOCNet validated by the shooting method, in the case of 1 agent.} Our solution matches the solution provided by shooting method, which means the result is a solution to the Hamiltonian ODE with $H$ provided by \eqref{eq:hamiltonian2}.}
    \label{fig:1d_valid}
\end{figure}

We further validate our training framework on this low-dimensional problem by benchmarking against the shooting method. It can be seen in~\Cref{fig:1d_valid} that our solution matches the solution provided by the shooting method, which means that the result is a solution to the Hamiltonian ODE with $H$ provided by \eqref{eq:hamiltonian2}.

\begin{figure}[h]
    \centering
    \includegraphics[width=0.8\textwidth]{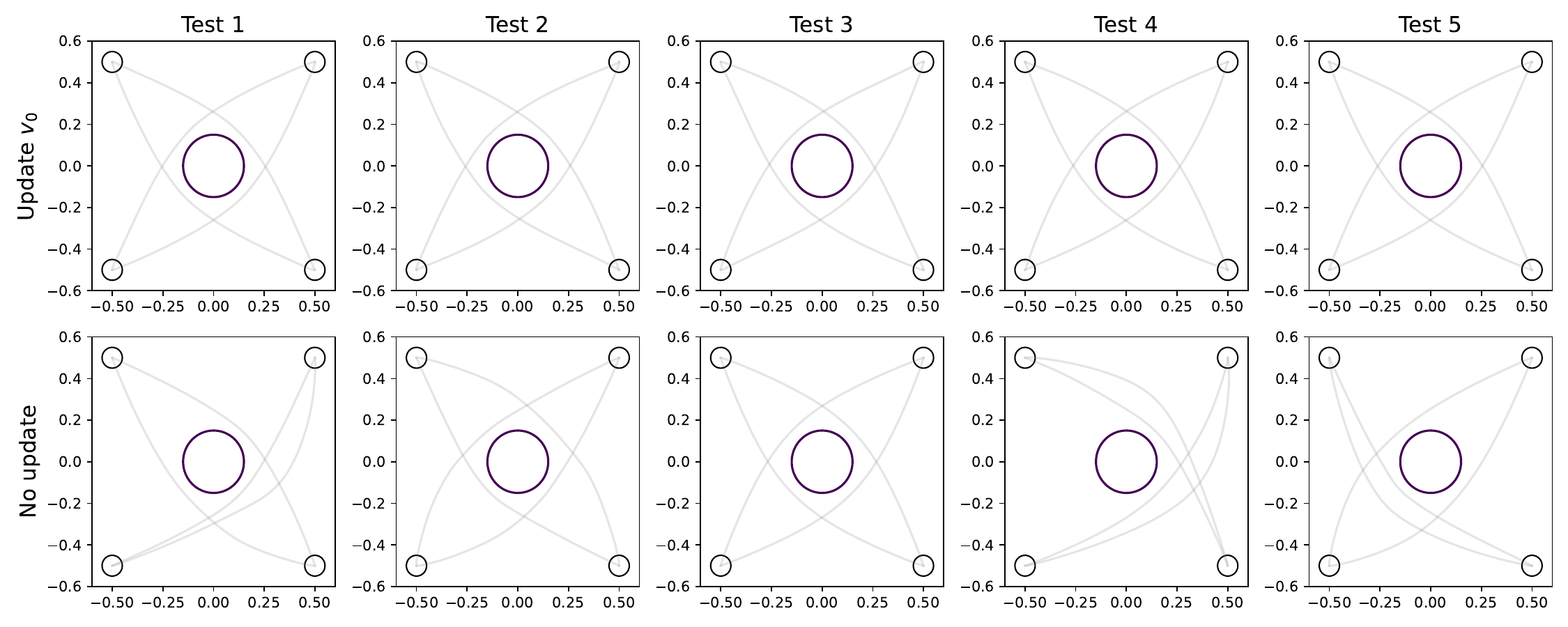}
    \caption{\textbf{Comparison of planned trajectory with or without the warm-up, at 5 different initializations.} The figure shows that our algorithm is more robust to random initialization and converges more often to global minima when the warm-up scheduler is applied.}
    \label{fig:4d_traj}
\end{figure}

\begin{figure}[h]
    \centering
    \includegraphics[width=0.8\textwidth]{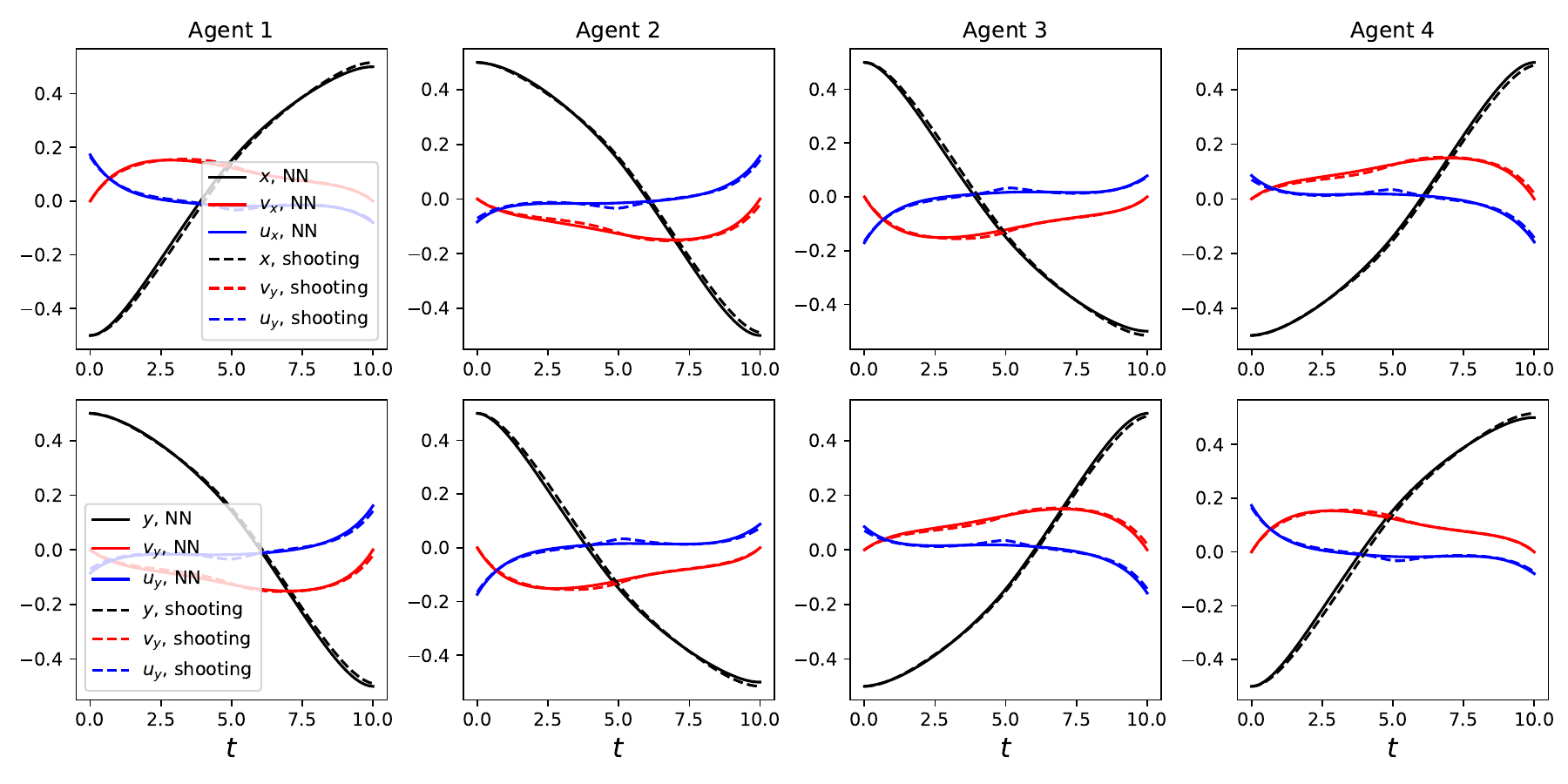}
    \caption{\textbf{The solution generated by our algorithm validated by the shooting method, in the case of 4 agents.} Our solution matches the solution provided by shooting method, which means the result is a solution to the Hamiltonian ODE.}
    \label{fig:4d_valid}
\end{figure}

\textbf{Four agents}
Assume $M = 4$, $L(\bx, \bu) = \frac{1}{2} \sum_{i=1}^4 \left(\|\bu_i\|^2 +\|\bv_i\|^2\right)$, $k = 0$. The initial positions of the drones are near the boundary of the room, and the terminal positions are the opposite locations, i.e., we set $\bw_T = -\bw_0$. It is also assumed that the agents start and terminate with zero velocity, as in the previous example. This problem contains several asymmetric local optima and two circular symmetric global optima. We found that with the training warm-up, our method converges to the global optima in all 5 independent simulations. However, if the training warm-up is not applied, 3 out of 5 cases TSympOCNet converge to the local optima, as shown in~\Cref{fig:4d_traj}. The solution is further validated against the shooting method, which exhibits great alignment, demonstrating that it is a solution to the Hamiltonian ODE, as shown in~\Cref{fig:4d_valid}. Note that we found that the shooting method does not converge well for the $2Md_x = 32$ dimensional Hamiltonian ODE. So we first use the shooting method to find the solution for one agent without considering the mutual collision avoidance requirement characterized by $h_2$. Then we rotate that solution by $90^\circ, 180^\circ, 270^\circ$, respectively, to obtain the solution for all 4 agents.

\subsection{High dimensional problem with Newtonian dynamics}\label{sec:highd}
\begin{figure}[h]
    \centering
    \includegraphics[width=0.8\textwidth]{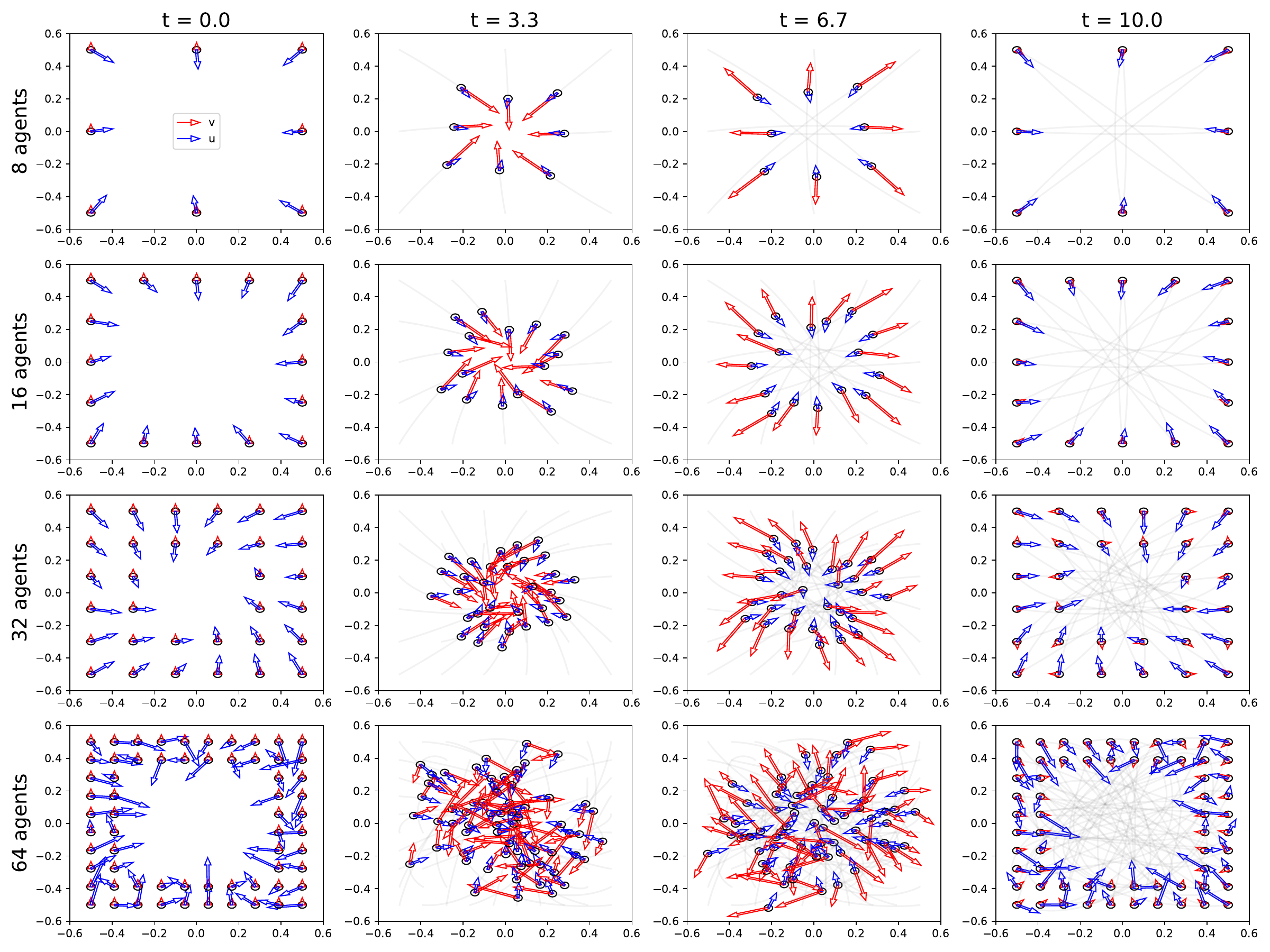}
    \caption{\textbf{The planned path and control for $M = 8, 16, 32, 64$ agents and time step $t = 0, 3.3, 6.7, 10$.} Animations for these 4 scenarios can be found in \href{https://github.com/zzhang222/TSympOCNet/tree/main/4.2a}{https://github.com/zzhang222/TSympOCNet/tree/main/4.2a}.}
    \label{fig:high_dim_traj}
\end{figure}

We consider the path planning problem of agents in a room of size $[-0.5, 0.5] \times [-0.5, 0.5]$. Collisions among the agents and collisions between agents and the room walls are prevented. The function $D$ in \eqref{eqt:def_d} is defined by:
\begin{equation*}
D_{bd}(\bw) = (w_1 + 0.5, 0.5-w_1, w_2 + 0.5, 0.5-w_2) \quad \forall \bw=(w_1,w_2)\in\R^2.
\end{equation*} 
In this experiment, we set $C_d = 0.02$. We assume $L(\bx, \bu) = \frac{1}{2}\|\bu\|^2$, $k = 0$.
The drones' initial positions lie near the room boundary, and their terminal positions are opposite, denoted as $\bw_T = -\bw_0$. This constitutes a high-dimensional example (state space dimension $n=4M$, varying from $32$ to $256$ in our experiments).
We run 10 repeated experiments with number of agents $M = 8,16,32,64$ to test the robustness of TSympOCNet. We further compare the solution obtained from TSympOCNet with solution obtained from PINN. For both neural network architectures, we use the same training algorithm in~\Cref{sec:training}, both with training warmup scheduler. The only difference between two approaches lies in the parameterization of $(\bx_\theta(s), \bp_\theta(s))$. TSympOCNet represents it by composition of latent LQR solution and TL-SympNet, while PINN sets $(\bx_\theta(s), \bp_\theta(s))$ as a fully-connected neural network on $s$.
In each trial, we document the training time in seconds, the running cost $L(\bx, \bu)$ and the constraint violation metric
 \begin{equation}
     \mathcal{D} = \min_s \min_{1\leq i < j \leq M} ||\bw_i(s) - \bw_j(s)|| - 2C_d,
 \end{equation}
where $\bw_i(s)$ denotes the center position for the $i$-th drone at time $s$. The mean and standard deviation of these three metrics are reported in~\Cref{tab:highd_stats}. It can be seen that TSympOCNet consistently outperforms the vanilla PINN in terms of generating a trajectory with lower running cost and lower constraint violation. It is worth remarking that in the case of $M = 8, 16$ agents, the constraint is exactly satisfied, which means that no collision is going to happen. However, in the case of 64 agents, the constraint is slightly violated in most of the solutions provided by TSympOCNet. To provide a feasible solution, one needs to rescale the problem and set $C_d = 0.16$. On the other hand, collision always occurs in PINN solutions when $M = 64$. The solution trajectory by TSympOCNet in one trial is shown in~\Cref{fig:high_dim_traj}. 

\begin{figure}[h]
    \centering
    \includegraphics[width=0.8\textwidth]{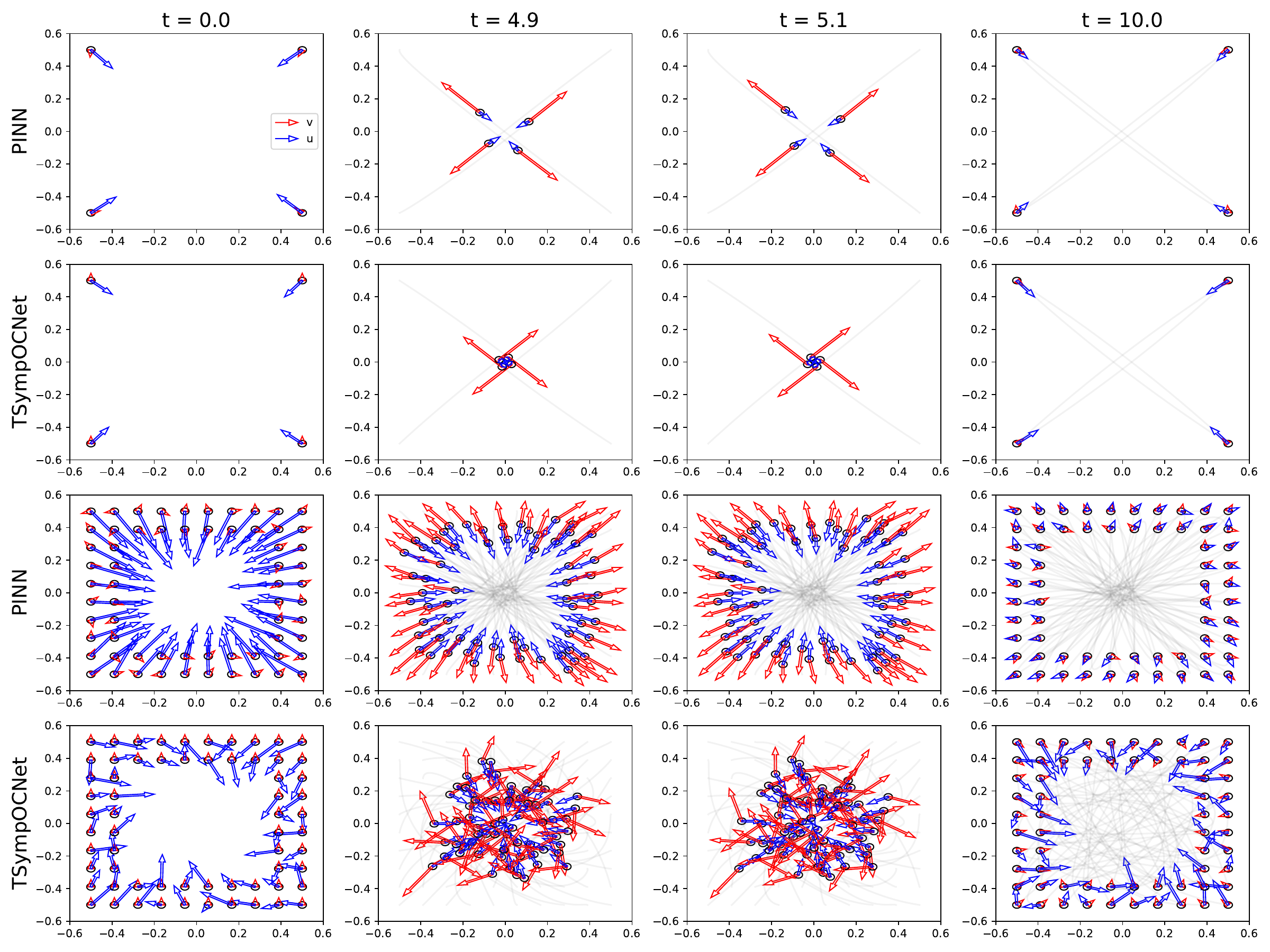}
    \caption{\textbf{The planned path and control for $M = 4$ and $64$ agents using PINN and TSympOCNet.} TSympOCNet produces a solution with lower cost in lower-dimensional problems and more feasible solution in higher dimensions. Animations for the solutions provided by TSympOCNet and PINN can be found in \href{https://github.com/zzhang222/TSympOCNet/tree/main/4.2b}{https://github.com/zzhang222/TSympOCNet/tree/main/4.2b}.}
    \label{fig:pinn_vs_sympocnet}
\end{figure}

We dive deeper into the solutions provided by PINN and TSympOCNet, focusing on cases where $M = 4$ and $M = 64$, as illustrated in~\Cref{fig:pinn_vs_sympocnet}. When $M = 4$, TSympOCNet's planned path exhibits a distinct behavior: the agents initially converge towards the origin, then just before encountering one another, initiating a synchronized rotation before dispersing. Conversely, PINN's path lacks this rotational symmetry—One pair of agents accelerates, passing ahead, while the second pair decelerates, awaiting the first pair's passage before proceeding. The first solution is more optimal in terms of running cost. In scenarios where $M = 64$, TSympOCNet's planned path demonstrates fewer constraint violations compared to PINN. The computation cost does not scale with the number of agents, since we use a fixed number of layers and fixed width for neural networks, in all the experiments of this subsection.

\begin{table}[h]
\tiny
\centering
\begin{tabular}{|cc|c|c|c|c|}
\hline
\multicolumn{2}{|c|}{\# agents}                                          & 8 & 16 & 32 & 64 \\ \hline
\multicolumn{1}{|c|}{\multirow{3}{*}{PINN}}       & Runtime (s)  &  $338\pm2$  & $338\pm2$   &  $354\pm2$  &  $368\pm2$  \\ \cline{2-6} 
\multicolumn{1}{|c|}{}                            & $L(\bx, \bu)$       &  $0.080\pm0.007$ & $0.175\pm0.013$   &  $0.346\pm0.020$  &   $4.516\pm1.091$ \\ \cline{2-6} 
\multicolumn{1}{|c|}{}                            & $\mathcal{D}$ &  $(-1 \pm 2)\times 10^{-4}$ &  $(-4 \pm 6)\times 10^{-4}$  &  $(-2 \pm 3)\times 10^{-3}$  &  $(- 39 \pm 1)\times 10^{-3}$  \\ \hline
\multicolumn{1}{|c|}{\multirow{3}{*}{TSympOCNet}} & Runtime (s)  & $465\pm2$  & $504\pm2$   &  $588\pm67$  &  $529\pm51$  \\ \cline{2-6} 
\multicolumn{1}{|c|}{}                            & $L(\bx, \bu)$     &  $0.072\pm0.0003$ & $0.152\pm0.006$   &  $0.304\pm0.049$  &   $1.125\pm0.118$ \\ \cline{2-6} 
\multicolumn{1}{|c|}{}                            & $\mathcal{D}$ &  0 &  0  &  $(-2 \pm 4)\times 10^{-4}$  &  $(-8 \pm 1)\times 10^{-3}$  \\ \hline
\end{tabular}
\caption{\textbf{Comparison of vannila PINN and TSympOCNet on problems of dimension $8,16,32,64$.} Under the provided hyperparameter setup, the solution obtained by TSympOCNet is closer to the global optima compared to vannila PINN. The collision-avoidance constraint is also satisfied better.}
\label{tab:highd_stats}
\end{table}

\subsection{Planning in non-convex environment} \label{sec:nonconvex}
\begin{figure}[h]
    \centering
    \includegraphics[width=0.8\textwidth]{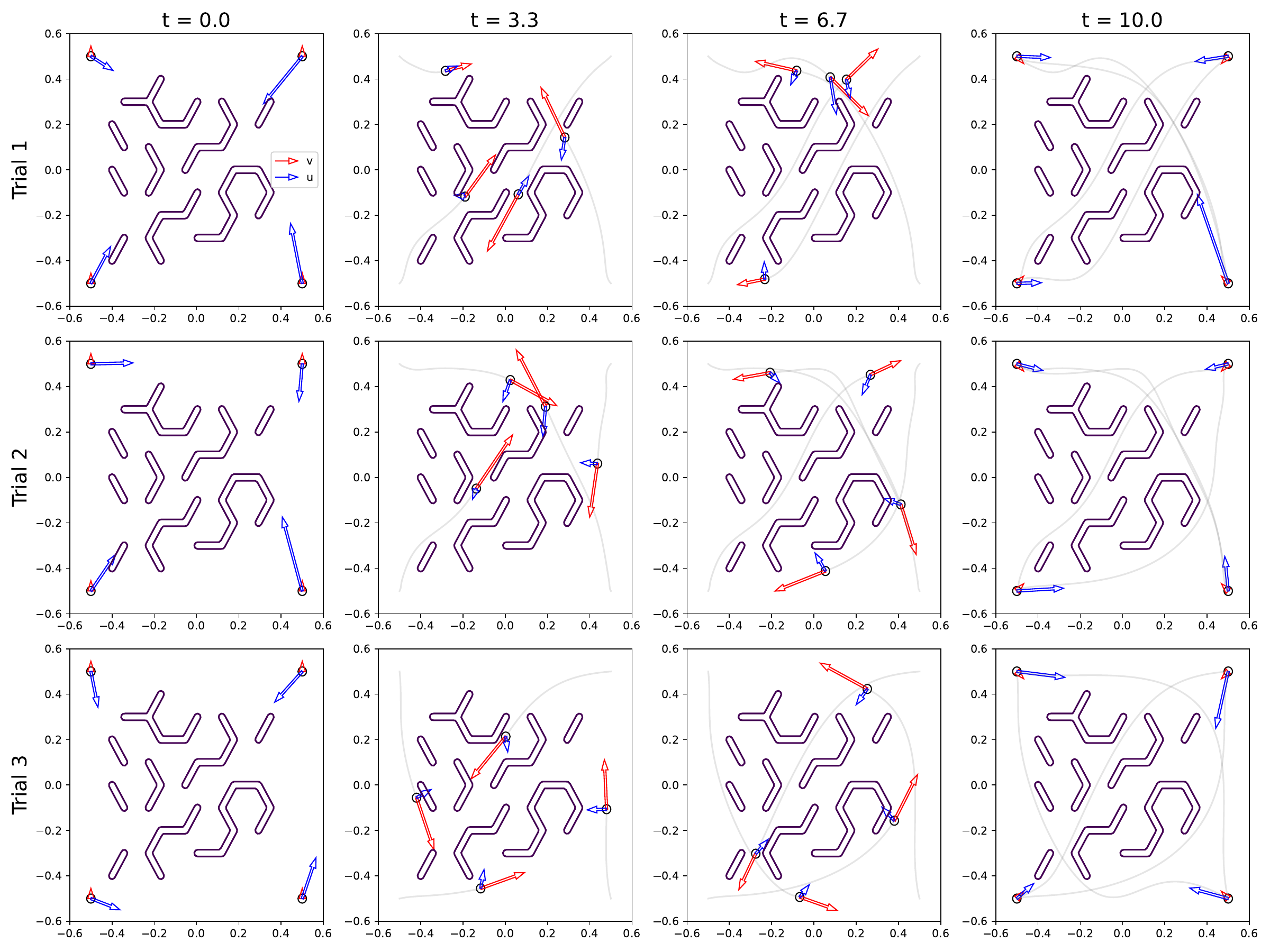}
    \caption{\textbf{The planned path and control in the non-convex maze.} We run three independent simulations and obtained three different solutions. The running cost $L(\bx, \bu) = 0.104, 0.151, 0.193$ in trial $1, 2, 3$. In all cases, no constraint violations are observed. TSympOCNet may converge to suboptimal solutions in this numerical example. Animations for these 3 trials can be found in \href{https://github.com/zzhang222/TSympOCNet/tree/main/4.3}{https://github.com/zzhang222/TSympOCNet/tree/main/4.3}.}
    \label{fig:maze}
\end{figure}

In this experiment, we set $C_d = 0.02$. We assume $L(\bx, \bu) = \frac{1}{2}\|\bu\|^2$, $k = 0$. We consider the path planning problem of agents in a room of size $[-0.5, 0.5]\times[-0.5,0.5]$. Apart from collisions among the agents and collisions between agent and the room boundary described in the previous section, we also aim to prevent the collision between agents and several walls inside the room, which created a maze. The shape of the maze can be seen in~\Cref{fig:maze}.
\begin{equation*}
\begin{split}
D(\bw) &= (D_{bd}(\bw), D_{maze}(\bw)) \\
D_{maze}^{(j)}(\bw) &=  \min_{\by\in l_j}\|\bw-\by\| - (C_o+C_d) \quad \forall \bw\in\R^2.
\end{split}
\end{equation*} 

The expression $\min_{\by\in l_j}\|\bw-\by\|$ in $D_{maze}$ calculates the distance between point $\bw$ and the line segment $l_j$. Considering that obstacle $E_j$ encompasses all points within a distance of $C_o$ from line segment $l_j$, a drone collides with $E_j$ if and only if the distance between its center and $l_j$ is less than $C_o + C_d$. Thus, the function $D_{maze}^{(j)}$ defined in this manner imposes a constraint that prevents drone collisions with the $j$-th obstacle $E_j$. 

We ran three independent simulations and observed three different solutions, as shown in~\Cref{fig:maze}. The running cost $L(\bx, \bu) = 0.104, 0.151, 0.193$ in trial $1, 2, 3$. In all cases, no constraint violations are observed. It is worth noting that TSympOCNet may converge to suboptimal solutions in this numerical example.
\subsection{Collision and obstacle avoidance in three-dimensional space} \label{sec:numerical_3d} 
In this experiment, we consider a path planning problem where each agent operates in dimension 3. 
This example is inspired by the swarm path planning problem studied in~\cite{onken2021neural, meng2022sympocnet}. 
Here we consider in the Newtonian dynamics, 
that is we further control the acceleration of each agent. 
Other setup are similar as in~\cite{onken2021neural,meng2022sympocnet}. 
In more detailed, we consider $M= 100$ drones with radius $0.2$, 
which leads to the problem in a dimension $6M=600$. We assume $L(\bx, \bu) = \frac{1}{2}\|\bu\|^2$, $k = 0$. 
Two three-dimensional rectangular obstacles are placed between the initial and the final positions. 
This is represented by the constraint function $h_1$. 
In particular, denote $[C^j_{11}, C^j_{1,2}] \times [C^j_{21},C^j_{22}] \times [C^j_{31}, C^j_{32}]$ the $j$th rectangular obstacle, the function $D_j$ in~\eqref{eqt:def_d} is defined by 
\begin{equation}
D_j(\bw) = \max_{i=1,2,3}\{C^j_{i1} - C_d - x_i, x_i - C^j_{i2} -C_d  \} , \forall \ \bw=(w_1,w_2,w_3)\in\R^3 \ .
\end{equation}
Notice that in this setup, the condition $D_j (\bw)\geq 0$ also ensures that the drone are at a distance $C_d$ away from each other, that is the collision avoidance is also represented in $h_1$. 
In numerical implementation, we set 
\begin{equation}
\begin{aligned}
& (C^1_{11},C^1_{12},C^1_{21},C^1_{22},C^1_{31},C^1_{32}) = (-1.8,1.8,-0.3,0.3,0.2,6.8) \\ 
& (C^2_{11},C^2_{12},C^2_{21},C^2_{22},C^2_{31},C^2_{32}) = (2.2,3.8,-0.8,0.8,0.2,3.8)
\end{aligned}
\end{equation}

\begin{figure}[htbp]
    \centering
    \begin{subfigure}{0.32\textwidth}
        \centering \includegraphics[width=\textwidth]{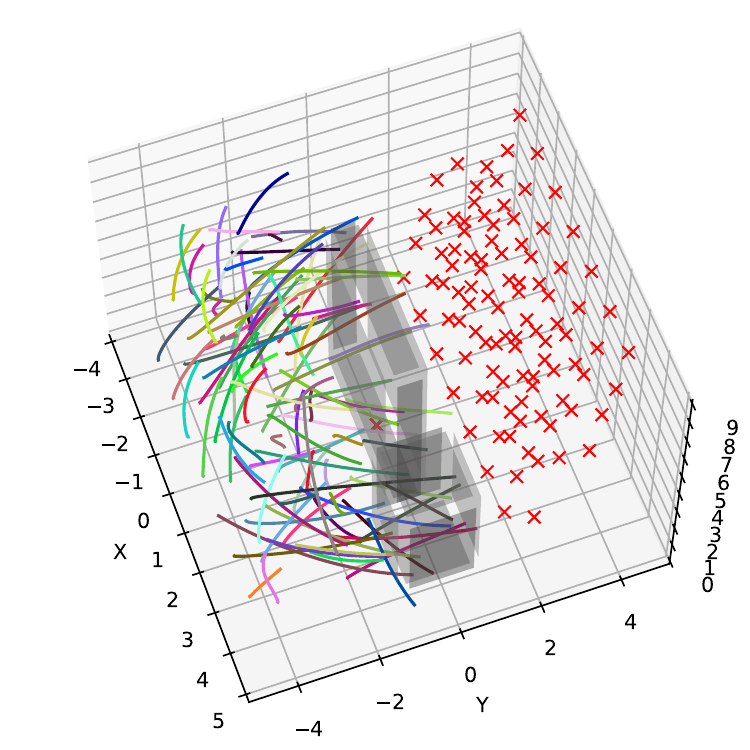}
        \caption{$t=3.3$}
    \end{subfigure}
    \begin{subfigure}{0.32\textwidth}
        \centering \includegraphics[width=\textwidth]{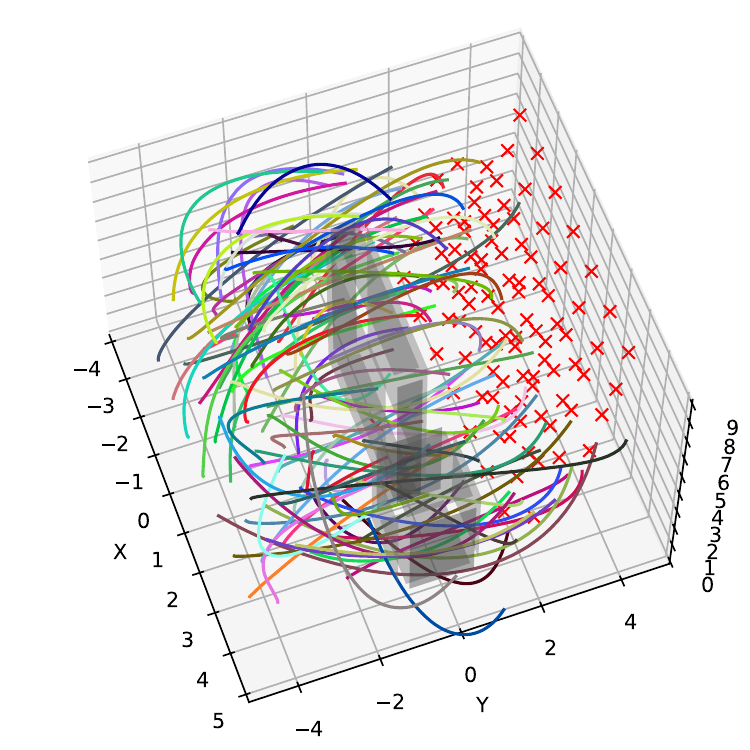}
        \caption{$t=6.7$}
    \end{subfigure}
    \begin{subfigure}{0.32\textwidth}
        \centering \includegraphics[width=\textwidth]{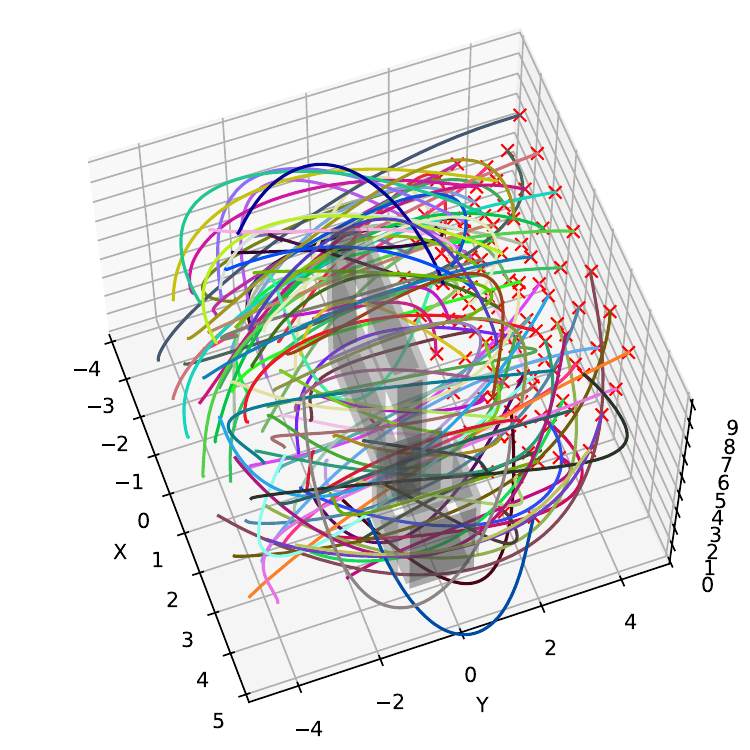}
        \caption{$t=10$}
    \end{subfigure}
    \caption{\textbf{Path planning of 100 drones in a 3d space.} We plot the predicted positions of 100 drones at time $t=3.3,6.7,10$. The paths from the initial positions to the current positions of all drones are plotted as colored lines. The destinations are marked as red crosses. The drones reach their destinations without any collision. The animation can be found in \href{https://github.com/zzhang222/TSympOCNet/blob/main/4.4/swarm.gif}{https://github.com/zzhang222/TSympOCNet/blob/main/4.4}.}
    \label{fig:figure_3d}
\end{figure}

The results are plotted in~\Cref{fig:figure_3d}.
\section{Summary}\label{sec:conclusion}
We introduce TSympOCNet, an extension of SympOCNet for tackling high-dimensional optimal control problems with state constraints and more general dynamics. Applying TSympOCNet to multi-agent simultaneous path planning tasks with obstacle avoidance demonstrates its efficacy in solving high-dimensional problems in hundreds of dimensions. These findings highlight TSympOCNet's potential for real-time solutions in high-dimensional optimal control problems. In future research, we aim to explore extending the current framework to handle feedback control problems by solving the latent LQR problem via Riccati solvers.

\section*{Acknowledgement}
The simulations were run on H100 GPU by NVIDIA. The discussion with Prof. Yeonjong Shin and Mr. Qian Zhang is greatly appreciated.

\newcommand{\etalchar}[1]{$^{#1}$}

\end{document}